\documentclass[12pt]{article}
\usepackage{amssymb,latexsym,epsfig}
\font\of=msbm10 scaled 1200
\def\R{\mbox{\of R}}
\def\C{\mbox{\of C}}
\def\Z{\mbox{\of Z}}
\def\N{\mbox{\of N}}

\newtheorem{definition}{Definition}

\newtheorem{theorem}{Theorem}%[section]
\newtheorem{lemma}{Lemma}
\newtheorem{proposition}{Proposition}

\newtheorem{corollary}{Corollary}
\title{The displacement map associated to polynomial unfoldings
of planar Hamiltonian vector fields}

 \author{Lubomir Gavrilov \\
 \normalsize \it Laboratoire Emile Picard, CNRS UMR 5580,  Universit\'e
 Paul Sabatier\\
 \normalsize \it 118, route de Narbonne, 31062 Toulouse Cedex, France \\
 \\ Iliya D. Iliev\\
 \normalsize \it Institute of Mathematics, Bulgarian Academy of Sciences\\
 \normalsize \it Acad. G. Bonchev Str, Bldg 8,
 1113 Sofia, Bulgaria } \date{September 10, 2004}

\begin{document}
\maketitle
%\table
\begin{abstract}{We study the displacement map associated to
small one-parameter polynomial unfoldings of polynomial
Hamiltonian vector fields on the plane. Its leading term,
the generating function $M(t)$, has an analytic continuation
in the complex plane and the real zeroes of $M(t)$ correspond
to the limit cycles bifurcating from the periodic orbits of
the Hamiltonian flow. We give a geometric description of the
monodromy group of $M(t)$ and use it to formulate sufficient
conditions for $M(t)$ to satisfy a differential equation of
Fuchs or Picard-Fuchs type. As
examples, we consider in more detail the Hamiltonian vector fields
$\dot{z}=i\bar{z}-i(z+\bar{z})^3$ and $\dot{z}=iz+\bar{z}^2$,
possessing a rotational symmetry of order two and three, respectively.
In both cases $M(t)$ satisfies a Fuchs-type equation but
in the first example $M(t)$ is always an Abelian integral
(that is to say, the corresponding equation is of Picard-Fuchs type)
 while
in the second one this is not necessarily true. We derive an explicit
formula of $M(t)$ and estimate the number of its real zeroes.}
\end{abstract}

\section{Introduction}

Consider a perturbed planar Hamiltonian vector field
$$ \left\{ \begin{array}{l}
\dot{x}  =   H_y(x,y) + \varepsilon P(x,y,\varepsilon ),\\
\dot{y}  =  -H_x(x,y) + \varepsilon Q(x,y,\varepsilon ).
\end{array}
\right.
\eqno{(1_\varepsilon )}
$$
We suppose that $H,P,Q$ are real polynomials in $x,y$ and moreover,
$P,Q$ depend analytically on a small real parameter $\varepsilon$.
%%%%%%%%%%%%%
Assume that for a certain open interval $\Sigma\subset\R$, the level sets of
the Hamiltonian $\{H=t\}$, $t\in\Sigma$, contain a continuous in $t$ family
of ovals $\cal A$. (An oval is a smooth simple closed curve which is free of
critical points of $H$). Such a family is called a {\it period annulus} of
the unperturbed system $(1_0)$. Typically, the endpoints of $\Sigma$ are
critical levels of the Hamiltonian function that correspond to centers,
saddle-loops or infinity. The limit cycles (that is, the isolated periodic
trajectories) of $(1_\varepsilon )$ which tend
to ovals from $\cal A$ as $\varepsilon\to 0$ correspond to the zeros of the
{\it displacement map} ${\cal P}_{\varepsilon }(t)-t$, where the
first return map ${\cal P}_{\varepsilon }(t)$ is defined on Fig. \ref{return}.
More explicitly, take a segment $\sigma$ which is transversal to the family
of ovals $\cal A$ and parameterize it by using the Hamiltonian value $t$.
For small $\varepsilon$, $\sigma$ remains transversal to the flow of
$(1_\varepsilon)$, too. Take a point $S\in\sigma$ and let $t=H(S)$.
The trajectory of $(1_\varepsilon)$ through $S$, after making one round,
will intersect $\sigma$ again at some point $S_1$ and the first return
map ${\cal P}_{\varepsilon }(t)$ is then defined by $t\rightarrow H(S_1)$.
%%%%%%%%%%%%%

Fixing a period annulus $\cal A$
of $(1_0)$ and taking a nonintegrable deformation $(1_\varepsilon)$,
then the related displacement map is defined in the corresponding open
interval $\Sigma\subset\R$ and there is a natural number $k$ so that
$$
{\cal P}_{\varepsilon }(t) - t = M(t) \varepsilon ^k +...,\quad t\in\Sigma.
\eqno{(2_k)}
$$
The limit cycles of $(1_\varepsilon )$ which tend to periodic
orbits from $\cal A$ as $\varepsilon\to 0$ correspond therefore
to the zeros of the {\it generating function} $M(t)$ in $\Sigma$.

The goal of the paper is to study the analytic continuation of the
generating function $M(t)$ in a complex domain. We give a geometric
description of the monodromy group of $M(t)$ (Theorem \ref{th1}) from
which we deduce sufficient conditions for $M(t)$ to satisfy a differential
equation of Fuchs or Picard-Fuchs type (Theorem \ref{main}).

Recall that a Fuchsian equation is said to be of Picard-Fuchs
type, provided that it possesses a fundamental set of solutions
which are Abelian integrals (depending on a parameter). In the
present paper by an Abelian integral we mean a function of the
form
\begin{equation}\label{ai} I(t)= \int_{\delta (t)} \omega
\end{equation}
where
\begin{itemize}
    \item  $\omega$ is a rational one-form in $\C^2$;
    \item there exists a bivariate polynomial $f: \C^2 \rightarrow
    \C$ such that $\delta (t) \subset f^{-1}(t)$, where $\{\delta
    (t)\}$ is a  family of closed loops, depending continuously on
    the complex parameter $t$.
\end{itemize}
It is supposed that $t$ belongs to some simply connected open
subset of $\C$ and $\delta (t)$ avoids the possible singularities
of the one-form $\omega$ restricted to the level sets $f^{-1}(t)$.
Under these conditions  $I(t)$ satisfies a linear differential
equation of Fuchs, and hence of Picard-Fuchs type.

It is well known that for a generic perturbation in $(1_\varepsilon)$
one has $k=1$ in $(2_k)$ and moreover,
$$
M(t)= \int_{\delta (t)} Q(x,y,0)dx -P(x,y,0)dy, \quad t\in\Sigma
$$
is then an Abelian integral \cite{pont}. Here $\delta (t)\subset
\R^2$, ${\cal A}=\{\delta(t)\}$, $t\in\Sigma$, is the continuous
family of ovals defined by the polynomial $H(x,y)$ and the
monodromy of $M(t)$ is deduced from the monodromy of $\delta (t)$
in a complex domain. More precisely, let $\Delta $ be the finite
set of atypical values of $H: \C^2\rightarrow \C$. The homology
bundle associated to the polynomial fibration
$$
\C^2\setminus H^{-1}(\Delta ) \stackrel{H}{\rightarrow } \C\setminus \Delta
$$
has a canonical connection. The monodromy group of the Abelian
integral $M(t)$ is then the monodromy group of the connection (or
a subgroup of it). It is clear that $M(t)$ depends on the homology
class of $\delta (t)$ in $H_1(\Gamma_t ,\Z)$ where $\Gamma_t $ is
the algebraic curve $\{(x,y)\in \C^2: H(x,y)=t \}$.

On the other hand, there are perturbations $(1_\varepsilon)$ with
$k>1$ in $(2_k)$. This happens when the perturbation is so chosen
that the first several coefficients in the expansion of the
displacement map, among them the function $M(t)$ given by the
above explicit integral, are identically zero in $\Sigma$. One
needs to consider such perturbations in order to set a proper
bound on the number of bifurcating limit cycles e.g. when the
Hamiltonian possesses symmetry or the degree of the perturbation
is greater than the degree of the original system. Therefore, the
case when $k>1$ is the more interesting one, at least what
concerns the infinitesimal Hilbert's 16th problem which is to find
the maximal number of limit cycles in $(1_\varepsilon)$, in terms
of the degrees of $H,P,Q$ only. In this case the generating
function $M(t)$ can have more zeroes in $\Sigma$, and respectively
the perturbations with $k>1$ can produce in general more limit
cycles than the ones with $k=1$ (see e.g. \cite{advde},
\cite{ili98}, \cite{gi00} for examples). Moreover, this case is
more difficult because the generating function is not necessarily
an Abelian integral and even the calculation of $M(t)$ itself is a
challenging problem. It turns out that in general (when $k>1$),
the generating function $M(t)$ depends on the {\em free homotopy
class} of the closed loop $\delta (t) \subset \Gamma _t$
(Proposition \ref{continuation}). The homology group $H_1(\Gamma_t
,\Z)$ must be replaced in this case by another Abelian group
$H_1^\delta(\Gamma_t ,\Z) $ which we define in section
\ref{universal}. Although there is a canonical homomorphism
$$
H_1^\delta(\Gamma_t ,\Z)  \rightarrow H_1(\Gamma_t ,\Z)
$$
it is neither surjective, nor injective in general. The bundle
associated to $H_1^\delta(\Gamma_t ,\Z) $ has a canonical
connection too and this is the appropriate framework for the study
of $M(t)$. This construction might be of independent interest in
the topological study of polynomial fibrations.

To illustrate our results we consider in full details two examples
$$
H_{A_3}= \frac{y^2}{2} + \frac{(x^2-1)^2}{4}\;\;
\mbox{   and   }\;\;
H_{D_4}= x[y^2-(x-3)^2],
$$
that are known as the eight-loop Hamiltonian and the Hamiltonian
triangle. Note that $H_{A_3}$ and $H_{D_4}$ are deformations of
the isolated singularities of type $A_3$ and $D_4$ respectively,
chosen to possess a rotational symmetry of order 2 and 3. We
explain first how Theorem \ref{main} applies to these cases. In
the $A_3$ case the differential equation satisfied by the
generating function $M(t)$ is of Picard-Fuchs type. This means
that $M(t)$ is always an Abelian integral, as conjectured earlier
by the second author, see \cite{jmp02}. On the other hand, in the
$D_4$ case the equation is of Fuchs type and has a solution which
is not a linear combination of Abelian integrals of the form
(\ref{ai}), with $f= H_{D_4}$. The reason is that the generating
function $M(t)$ has a term $(\log(t))^2$ in its asymptotic
expansion. Equivalently, the monodromy group of the associated
connection contains an element of the form
$$
 \pmatrix{
1 & 1& 0  \cr 0 & 1& 1 \cr 0 & 0& 1  }
$$
which could not happen if the associated equation were of
Picard-Fuchs type. Next, we provide an independent study of $M(t)$
based on a generalization of Fran\c{c}oise's algorithm
\cite{fra96}. It is assumed for simplicity that in
$(1_\varepsilon)$ the polynomials $P, Q$ do not depend on
$\varepsilon$. In the $A_3$ case, we derive explicit formulas for
$M(t)$ in terms of $k$ and the degree $n$ of the perturbation
(Theorem 3) and use them to estimate the number of bifurcating
limit cycles in $(1_\varepsilon)$ which tend to periodic orbits of
the Hamiltonian system (Theorems 4, 5, 6). Note that our argument
applies readily to the double-heteroclinic Hamiltonian
$H=\frac12y^2-\frac14(x^2-1)^2$ and to the global-center
Hamiltonian $H=\frac12y^2+\frac14(x^2+1)^2$ as well. What concerns
the Hamiltonian triangle, we give an explicit example of a
quadratic perturbation leading to a coefficient $M(t)$ at
$\varepsilon^3$ which is not an Abelian integral and derive the
third-order Fuchsian equation satisfied by $M(t)$. This part of
the paper uses only ``elementary" analysis and may be read
independently. We hope that the complexity of the combinatorics
involved will motivate the reader to study the rest of the paper.
This was the way we followed, when trying to understand the
controversial paper \cite{jmp02}
%%%%%%%%%%%%%%%%%
(its revised version is to appear in Bull. Sci. Math.).
%%%%%%%%%%%%%%%%%

The applications of Theorem \ref{main} which we present are by no means the
most general. On the contrary, these are the simplest examples in which it
 gives non-trivial answers. Theorem \ref{main} can be
 further generalized and a list of open questions is presented at the end of
 section \ref{mainresult}.

\section{Generating functions and limit cycles}

Assume that $f=f(x,y)$ is a real polynomial of degree at least 2 and
consider a polynomial foliation ${\cal F}_\varepsilon $ on the real plane
$\R^2$ defined by
\begin{equation}
\label{ee}
df - \varepsilon Q(x,y,\varepsilon )dx+\varepsilon P(x,y,\varepsilon )dy =0
\end{equation}
 where $P,Q$ are real polynomials in $x,y$ and analytic in $\varepsilon $, a
sufficiently small real parameter.
%%%%%%%%%%%%%%%
Note that (\ref{ee}) is just another form of the equation $(1_\varepsilon )$
with $H$ replaced by $f$.
%%%%%%%%%%%%%%%

Let $\delta (t)\subset \{(x,y)\in\R^2: f(x,y)=t \}$
be a continuous family of ovals defined on a maximal open interval
$\Sigma \subset \R$. We identify $\Sigma $ with a cross-section
$\Sigma \rightarrow  \R^2$ transversal to the ovals $\delta (t)$ from
the period annulus ${\cal A}= \cup_{t\in \Sigma } \delta (t)$.
For every compact sub-interval $K \subset \Sigma $, there exists
$\varepsilon _0 = \varepsilon _0(K)$ such that the first return map
${\cal P}_\varepsilon (t)$ associated to the period annulus ${\cal A}$
is well defined and analytic in
$$
\{ (t,\varepsilon ) \in \R^2: t\in K, |\varepsilon | < \varepsilon _0 \} \; .
$$
As the limit cycles of (\ref{ee}) intersecting $K$ correspond to the
isolated zeros of ${\cal P}_\varepsilon (t)-t$, we shall always suppose
that ${\cal P}_\varepsilon (t) \not \equiv t$.
Then there exists $k\in \N$ such that
\begin{equation}
\label{displ}
{\cal P}_\varepsilon (t)-t = M_k(t) \varepsilon ^k + O(\varepsilon ^{k+1})
\end{equation}
uniformly in $t$ on each compact sub-interval $K$ of $\Sigma $. Therefore the
number of the zeros of $M_k(t)$ on $\Sigma$ provides an upper
bound to the number of zeros of ${\cal P}_\varepsilon (t)-t$ on $\Sigma$ and
hence to the number of the corresponding limit cycles of
(\ref{ee}) which tend to $\cal A$ as $\varepsilon\to 0 $.
%%%%%%%%%%
Indeed, taking the right-hand side of (\ref{displ}) in the form
$\varepsilon ^k[M_k(t)  + O(\varepsilon)]$ and using the implicit function
theorem (respectively, the Weierstrass preparation theorem in the case of
multiple roots), we see that the displacement map and its first non-zero
coefficient $M_k(t)$ will have the same number of zeros in $\Sigma$ for
small $\varepsilon\neq 0$.
%%%%%%%%%%

\begin{definition}
\label{disp}
We call ${\cal P}_\varepsilon (t)-t$ the displacement map, and $M_k(t)$ the
$(k$-th$)$
 generating function, associated to the family of ovals $\delta(t)$
and to the unfolding ${\cal F}_\varepsilon $.
\end{definition}
{\bf Example.} If $f$  has $(\mbox{\rm deg}\; f-1)^2$
different critical points with different critical values, then
$M_k(t) = \int_{\delta (t)} \Omega _k$ where $\Omega _k$ is a polynomial
one-form in $x,y$. Therefore, the generating function $M_k(t)$ is an Abelian
integral. This easily follows from Fran\c{c}oise's  recursion formula
\cite{fra96} and the fact that if $ \int_{\delta (t)} \Omega  \equiv 0$
for a certain polynomial one-form $\Omega $, then $\Omega = dG+ gdf$ for
suitable polynomials $G,g$ \cite{ily78,gav98}. On the other hand, when
$f$ is non-generic (e.g. has ``symmetries"), this might not be true,
see the examples in section \ref{examples}.

\subsection{The monodromy group of the generating function}

For any non-constant complex polynomial $f(x,y)$  there exists a
finite set $\Delta \subset \C$ such that the fibration $\C^2
\stackrel{f}{\rightarrow } \C \setminus \Delta $ is locally
trivial. Let $t_0 \not \in \Delta$, $P_0\in f^{-1}(t_0)$ and
$\Sigma \subset \C^2 $ be a small complex disc centered at $P_0$
and transversal to $f^{-1}(t_0)\subset \C^2$. We will also suppose
that the fibers $f^{-1}(t)$ which intersect $\Sigma $ are regular,
hence $t=f(x,y)|_\Sigma $ is a coordinate on $\Sigma $.

To an unfolding ${\cal F}_\varepsilon $ of $df=0$ on the complex plane
$\C^2$ defined by (\ref{ee}), and to a closed loop
$$l_0: [0,1] \rightarrow  f^{-1}(t_0),\;\; l_0(0)=l_0(1)=P_0,$$
 we  associate a
holonomy map (return map,  Poincar\'e map in a complex domain)
$$
{\cal P}_{l_0,{\cal F}_\varepsilon } : \Sigma  \rightarrow \Sigma \;.
$$
\begin{figure}
\begin{center}
\epsfig{file=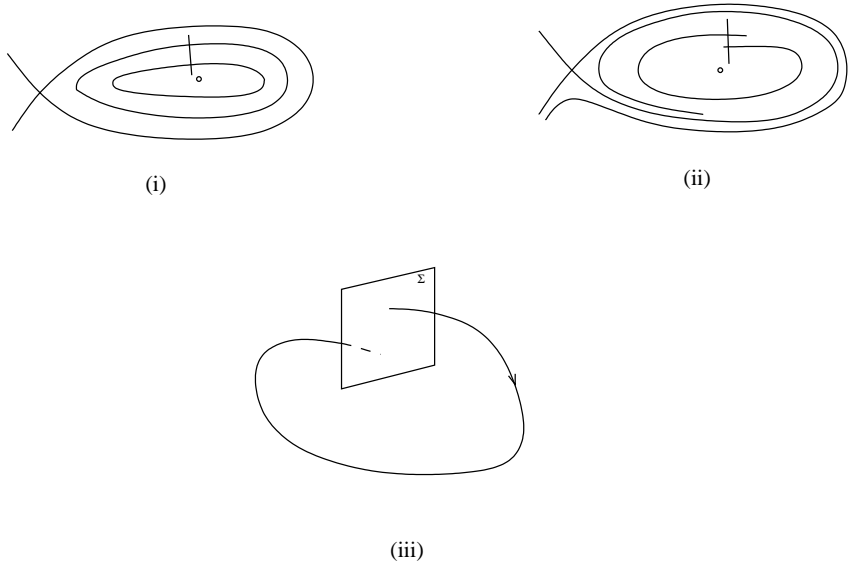}
\end{center}
\caption{The first return map and its complexification. }
\label{return}
\end{figure}
In the case when $l_0$ is an oval of the real polynomial $f$, it is just
the complexification of the analytic Poincar\'e map ${\cal P}_\varepsilon$
defined above, see Fig. \ref{return}. In  general, the definition of
${\cal P}_{l_0,{\cal F}_\varepsilon }$
is the following, see e.g. \cite{loray}.
Let ${\cal F}_0^{\perp}$ be a holomorphic foliation transversal to ${\cal
F}_0= \{ df=0 \}$ in some neighborhood of
$l_0$ (for instance, ${\cal F}_0^{\perp}=\{f_ydx - f_x dy=0\}$).
Then for $|\varepsilon |$ sufficiently small,
${\cal F}_0^{\perp}$ remains transversal to
${\cal F}_\varepsilon $. The holonomy map ${\cal P}_{l_0,{\cal
F}_\varepsilon }$ is a germ of a
biholomorphic map in a neighborhood of $P_0\in \Sigma $ which is obtained
by lifting the loop $l_0$
in the leaves of ${\cal F}_\varepsilon $ via  ${\cal F}_0^{\perp}$. Namely,
$Q= {\cal P}_{l_0,{\cal F}_\varepsilon }(P)$ if there exists a path
$\tilde{l_0}$ in a leaf of ${\cal F}_\varepsilon $ which connects
$P$ and $Q$, and which is a lift of the loop $l_0$ according
to ${\cal F}_0^{\perp}$. The holonomy map
${\cal P}_{l_0,{\cal F}_\varepsilon }$ does not depend on the choice of the
transversal foliation ${\cal F}_0^{\perp}$. If $l_0,l_1$ are two homotopic
loops with the same initial
point $P_0$, then ${\cal P}_{l_0,{\cal F}_\varepsilon } = {\cal
P}_{l_1,{\cal F}_\varepsilon }$.

%When there is no danger of confusion, we shall write from now on simply
%$$
%{\cal P}_{l_0,{\cal F}_\varepsilon }= {\cal P}_{l_0 }.
%$$
Let us fix the foliation ${\cal F}_\varepsilon$ and the loop
$l_0$.
 As before, if we suppose that ${\cal P}_{l_0,{\cal F}_\varepsilon }
\neq  id$, then there exists $k\in \N$ such that
$$
{\cal P}_{l_0,{\cal F}_\varepsilon }(t) = t+ \varepsilon ^k
M_{k}(l_0,{\cal F}_\varepsilon ,t)+...
$$
When there is no danger of confusion, we shall write simply
$$
M_{k}(l_0,{\cal F}_\varepsilon ,t)= M_k(t).
$$
The function $M_{k}$ is called {\em the generating function}
associated to the unfolding ${\cal F}_\varepsilon $ and to the
loop $l_0$. Note that the natural number $k$ as well as $M_k$ depend
on $l_0,{\cal F}_\varepsilon$ and $\Sigma$ in general. The following
observation is crucial for the rest of the paper.
\begin{proposition}
\label{sigma}
 The number $k$ and the generating function $M_k$ do not depend on
$\Sigma$. They depend on the foliation ${\cal F}_\varepsilon$ and
on the {\em free homotopy class} of the loop $l_0\subset
f^{-1}(t)$. The generating function $M_k(t)$ allows an analytic
continuation on the universal covering of $\C \setminus \Delta $,
where $\Delta$ is the set of atypical points of $f$.
\end{proposition}
The proof the proposition uses the following algebraic lemma.
\begin{lemma}
\label{continuation}
 Take $k\in\N$.  Let
$$
P_{\varepsilon}(t)= t + \sum_{i=k}^\infty \varepsilon^i p_i(t),
\quad p_k \neq 0, \quad
G_{\varepsilon}(t)= t + \sum_{i=1}^\infty \varepsilon^i g_i(t)
$$
be convergent power series of $(t,\varepsilon)$ in a suitable polydisc
centered at the origin in $\C^2$. If $\varepsilon$ is fixed and
sufficiently small, then $G_\varepsilon$ is a local automorphism and
$$
G_\varepsilon^{-1}\circ P_\varepsilon \circ G_\varepsilon (t)= t+
\sum_{i=k}^\infty \varepsilon^i \tilde{p}_i(t)
$$
where $ \tilde{p}_k(t) \equiv  p_k(t) $.
\end{lemma}
{\bf Proof of Lemma \ref{continuation}.} We have
\begin{eqnarray*}
  P_\varepsilon \circ G_\varepsilon (t) &=& G_\varepsilon (t) +
\sum_{i=k}^\infty \varepsilon^i p_i(G_\varepsilon (t))
  = G_\varepsilon (t) +  \varepsilon^k p_k(t) + {\rm O}(\varepsilon^{k+1})\\
 G_\varepsilon^{-1}(t)& = & t + \sum_{i=1}^\infty \varepsilon^i \tilde{g}_i(t)
\end{eqnarray*}
and therefore
\begin{eqnarray*}
 G_\varepsilon^{-1}\circ P_\varepsilon \circ G_\varepsilon (t)&=&
 G_\varepsilon (t) +  \varepsilon^k p_k(t) + {\rm O}(\varepsilon^{k+1}) \\
   & +& \sum_{i=1}^\infty \varepsilon^i \tilde{g}_i(G_\varepsilon (t) +
   \varepsilon^k p_k(t) + {\rm O}(\varepsilon^{k+1}))\\
   &=&G_\varepsilon (t) +  \varepsilon^k p_k(t)
   + \sum_{i=1}^\infty \varepsilon^i \tilde{g}_i(G_\varepsilon (t)
   ) + {\rm O}(\varepsilon^{k+1})\\
   &=& G_\varepsilon^{-1} \circ G_\varepsilon (t) + \varepsilon^k
   p_k(t) + {\rm O}(\varepsilon^{k+1}) \\
   & = & t + \varepsilon^k
   p_k(t) + {\rm O}(\varepsilon^{k+1}) .
\end{eqnarray*}
In the above computation ${\rm O}(\varepsilon^{k+1})$ denotes a power
series in $t,\varepsilon$ containing terms of degree at least
$k+1$ in $\varepsilon$.
The lemma is proved.\\

\vspace{2ex}
\noindent
{\bf Proof of  Proposition \ref{sigma}.} Let $\tilde{\Sigma}$ be
another transversal disc centered at $P_0$ and
$$
\tilde{{\cal P}}_{l_0,{\cal F}_\varepsilon }(t) : \tilde{\Sigma}
\rightarrow \tilde{\Sigma}
$$
the corresponding holonomy map. Then
$$
{\cal P}_{l_0,{\cal F}_\varepsilon }(t) = G_\varepsilon^{-1}\circ
\tilde{{\cal P}}_{l_0,{\cal F}_\varepsilon }(t) \circ
G_\varepsilon(t)
$$
where
$$
G_\varepsilon : \Sigma \rightarrow \tilde{\Sigma}
$$
is analytic and $G_0(t) \equiv t$. Lemma \ref{continuation} shows
that
$$
\tilde{{\cal P}}_{l_0,{\cal F}_\varepsilon }(t) = t +
\varepsilon^k
   M_k(t) + {\rm O}(\varepsilon^{k+1}),\quad M_k(t) \not\equiv 0
   $$
if and only if
$$
{\cal P}_{l_0,{\cal F}_\varepsilon }(t) = t + \varepsilon^k
   M_k(t) + {\rm O}(\varepsilon^{k+1}),\quad M_k(t) \not\equiv 0
   $$
As the holonomy map ${\cal P}_{l_0,{\cal F}_\varepsilon }(t)$
depends on the homotopy class of $l_0$ this holds true for $k$ and
$M_k$. In contrast to ${\cal P}_{l_0,{\cal F}_\varepsilon }$, the
generating function $M_k$ depends on the {\em free} homotopy class
of $l_0$. Indeed, let $\tilde{l}_0$ be a path in $f^{-1}(t_0)$
starting at $Q_0$ and terminating at $P_0$, and let
$\tilde{\Sigma}$ be a transversal disc centered at $Q_0$ with
corresponding holonomy map
$$
\tilde{{\cal P}}_{l_0,{\cal F}_\varepsilon }(t): \tilde{\Sigma}
\rightarrow \tilde{\Sigma} .
$$
Then we have
\begin{equation}\label{conjugation}
{\cal P}_{l_0,{\cal F}_\varepsilon }(t)= G_{\tilde{l}_0,{\cal
F}_\varepsilon } ^{-1}\circ \tilde{{\cal P}}_{l_0,{\cal
F}_\varepsilon }(t) \circ G_{\tilde{l}_0,{\cal F}_\varepsilon }(t)
\end{equation}
where
$$
G_{\tilde{l}_0,{\cal F}_\varepsilon } : \Sigma \rightarrow
\tilde{\Sigma}
$$
is analytic and $G_{\tilde{l}_0,{\cal F}_0 }(t) \equiv t$ (the
definition of $G_{\tilde{l}_0,{\cal F}_\varepsilon} $ is similar
to the definition of ${\cal P}_{l_0,{\cal F}_\varepsilon }(t)$).
Lemma \ref{continuation} shows that the generating function
$M_k(t)$
 does not depend on the
special choice of the initial point $P_0$. We conclude that it
depends only on the {\it free} homotopy class of the loop $l_0$.
Until now $M_k$ was defined only locally (on the transversal disc
$\Sigma$). As the fibration $\C^2 \setminus f^{-1}(\Delta )
\stackrel{f}{\rightarrow } \C \setminus \Delta $ is locally
trivial, then each closed loop $l_0 \in f^{-1}(t_0)$ defines
 a continuous family $l_0(t)$ of  closed loops on $f^{-1}(t)$, defined on
the universal covering space of $ \C \setminus \Delta $. Only the
free homotopy classes of the loops $l_0(t)$ are well defined and
to each $l_0(t)$ corresponds a holonomy map defined up to
conjugation, see (\ref{conjugation}). As this conjugation
preserves the number $k$ and the generating function $M_k(t)$ then
the latter allows an analytic continuation on the universal
covering of $\C \setminus \Delta $. Proposition \ref{sigma}
is proved. $\Box$

The monodromy group of  $M_k(t)$  is defined as follows. The
function $M_k(t)$ is multivalued on  $\C \setminus \Delta $. Let
us  consider all its possible determinations in a sufficiently
small neighborhood of $t=t_0$. All integer linear combinations of
such functions form a module over $\Z$ which we denote by ${\cal
M}_k( l_0,{\cal F}_\varepsilon )$. When there is no danger of
confusion we shall write simply
$$
{\cal M}_k( l_0,{\cal F}_\varepsilon ) = {\cal M}_k .
$$
The fundamental group $\pi _1(\C \setminus \Delta ,t_0) $ acts on
${\cal M}_k$ as follows. If $\gamma \in \pi _1(\C \setminus
\Delta,t_0)$ and $M \in {\cal M}_k$, let $\gamma _* M(t)$ be the
analytic continuation of $M(t)$ along $\gamma $. Then $\gamma _*$
is an automorphism of ${\cal M}_k$ and
$$(\gamma _{1}\circ \gamma _2)_* M= \gamma _{2*}(\gamma _{1*} M) \; .$$
\begin{definition}
\label{monodromy} The monodromy representation associated to the
generating function $M_k$ is  the group homomorphism
\begin{equation}
\label{repr} \pi _1(\C \setminus \Delta ,t_0) \rightarrow
Aut({\cal M}_k) \; .
\end{equation}
The group image of $\pi _1(\C \setminus \Delta ,t_0)$ under
$(\ref{repr})$ is called the monodromy group of $M_k$.
\end{definition}
In what follows we wish to clarify the case when the generating
function is (or is not) an Abelian integral. For this we need to
know the monodromy representation of $M_k$.

\subsection{The universal monodromy representation of the generating
function}
\label{universal}

Let $H$ be a  group and $S\subset H$ a set. We  construct an
 abelian group  $\hat{S}/[H,\hat{S}]$ associated to the pair $H,S$ as
follows. Let $\hat{S}$ be the group generated by the set
$$
\{ hsh^{-1}: h\in H\} ,
$$
that is to say, the least normal subgroup of $H$ containing $S$.
We denote by $[H,\hat{S}]$ the ``commutator" group generated by
$$
\{ hsh^{-1}s^{-1}: h\in H,\; s\in \hat{S} \} .
$$
Then $[H,\hat{S}]=[\hat{S},H]  $ is a normal subgroup of $\hat{S}$
and
$
 \hat{S}/[H,\hat{S}]
$ is an abelian group.  There is a canonical homomorphism
$$
\hat{S}/[H,\hat{S}] \rightarrow H/[H,H]
$$
which is not injective in general. Note that $\hat{S}=H$ implies
that $\hat{S}/[H,\hat{S}]= H/[H,H]$ is the abelianization of $H$.

We apply now the above construction to the case when $H=
\pi_1(\Gamma ,P_0)$ is the fundamental group of a connected
surface $\Gamma $ (not necessarily compact), $P_0\in \Gamma$. Let
$\pi_1(\Gamma )$ be the set of immersions of the circle into
$\Gamma $, up to homotopy equivalence (the set of free homotopy
classes of closed loops). Let $S\subset \pi_1(\Gamma)$ be a set
and $\hat{S} \subset \pi_1(\Gamma,P_0)$  be the pre-image of $S$
under the canonical projection
$$
 \pi_1(\Gamma,P_0) \rightarrow  \pi_1(\Gamma) \; .
$$
Then $\hat{S}$ is a normal subgroup of $\pi_1(\Gamma,P_0)$ and we
define
$$
H^S_1(\Gamma ,\Z)= \hat{S}/[\hat{S},\pi_1(\Gamma,P_0)] .
$$
In the case when $\hat{S} = \pi_1(\Gamma,P_0)$ we have
$H^S_1(\Gamma ,\Z) =  H_1(\Gamma ,\Z)$, the first homology group of
$\Gamma$. Let $\Psi$ be a diffeomorphism of $\Gamma $. It induces
a map
$$
\Psi_* : \pi_1(\Gamma ) \rightarrow \pi_1(\Gamma )
$$
and we suppose that $\Psi_*(S)=S$. Then it induces an automorphism (denoted
again by $\Psi*$)
$$
\Psi_* : H^S_1(\Gamma ,\Z) \rightarrow H^S_1(\Gamma ,\Z) .
$$
Note also that if $\Psi_0$ is a diffeomorphism isotopic to the identity,
then it induces the
identity automorphism.

Two closed loops $s_1,s_2 \in \hat{S}$ represent the same free
homotopy class if and only if $s_1 = h s_2 h^{-1}$ for some $h\in
\pi_1(\Gamma,P_0)$. It follows that to each free homotopy class of
closed loops  represented by an element of $\hat{S}$ there
corresponds a unique element of $H^S_1(\Gamma ,\Z)$.

Consider finally the locally trivial fibration
 $$\C^2 \setminus f^{-1}(\Delta ) \stackrel{f}{\rightarrow } \C\setminus
\Delta $$
 defined by the non-constant polynomial $f\in \C[x,y]$
and put $\Gamma = f^{-1}(t_0)$, $t_0 \not\in \Delta $.
Each loop $\gamma \in \pi_1(\C\setminus \Delta , t_0)$ induces a
diffeomorphism $\gamma _*$ of $\Gamma $,
defined up to an isotopy, and hence
a canonical group homomorphism
\begin{equation}
\label{class}
\pi_1(\C\setminus \Delta , t_0) \rightarrow \mbox{\rm Diff}\,(\Gamma
)/\mbox{\rm
Diff}_0(\Gamma )  .
\end{equation}
Here ${\rm Diff}\,(\Gamma )/\mbox{\rm
Diff}_0(\Gamma ) $ denotes the group of diffeomorphisms ${\rm
Diff}\,(\Gamma )$ of $\Gamma $, up to
diffeomorphisms ${\rm Diff}_0\,(\Gamma )$ isotopic to the identity
(the so called mapping class
group of $\Gamma $).
 The
homomorphism (\ref{class}) induces a homomorphism (group action on
$\pi_1(\Gamma )$)
\begin{equation}
\label{action1}
\pi_1(\C\setminus \Delta , t_0) \rightarrow \mbox{\it Perm}\,(\pi_1(\Gamma ))
\end{equation}
where $\mbox{\it Perm}\,(\pi_1(\Gamma ))$ is the group of permutations of
$\pi_1(\Gamma )$.

Let $l_0 \in \Gamma $ be a closed loop, and let $\hat{S} \subset
\pi_1(f^{-1}(t_0),P_0)$ be the
subgroup ``generated" by $l_0$. More precisely,
let $\bar{l_0}\in
\pi_1(f^{-1}(t_0))$ be the
free homotopy equivalence class represented by $l_0$. We denote by $S\subset
\pi_1(\Gamma )$ the orbit
$\pi_1(\C\setminus \Delta , t_0) \bar{l_0}$.
Let $\hat{S} \subset \pi_1(\Gamma ,P_0)$ be the subgroup generated by
the pre-image of the orbit
${\cal O}_{l_0}$ under the canonical map
$$
\pi_1(\Gamma ,P_0)  \rightarrow \pi_1(\Gamma ) \;
$$
and let us put
$$
H^{l_0}_1(\Gamma ,\Z)=  \hat{S}/[\pi_1(\Gamma ,P_0),\hat{S}].
$$
We obtain therefore the following
\begin{proposition}
\label{prdf} The group $H^{l_0}_1(f^{-1}(t_0) ,\Z)$ is abelian and
the canonical map
\begin{equation}
H^{l_0}_1(f^{-1}(t_0) ,\Z) \rightarrow  H_1(f^{-1}(t_0),\Z)
\label{canonical}
\end{equation}
is a homomorphism. The group action $(\ref{action1})$ of
$\pi_1(\C\setminus \Delta , t_0)$ on $\pi_1(f^{-1}(t_0))$
 induces a homomorphism
\begin{equation}
\label{monodromy1} \pi_1(\C\setminus \Delta , t_0) \rightarrow
Aut(H^{l_0}_1(f^{-1}(t_0) ,\Z))
\end{equation}
called the monodromy representation associated to the loop $l_0$.
\end{proposition}
The monodromy
group associated to $l_0$ is  the group image of $\pi _1(\C \setminus
\Delta ,t_0)$ under
the group homomorphism (\ref{monodromy1}).
\begin{theorem}
\label{th1} For every polynomial deformation ${\cal F}_\varepsilon
$ of the foliation $df=0$, and every closed loop $l_0 \subset
f^{-1}(t_0)$, the monodromy representation $(\ref{repr})$ of the
generating function $M_k$  is a sub-representation of the
monodromy representation dual to $(\ref{monodromy1})$.
\end{theorem}

The concrete meaning of the above theorem is as follows. There
exists a canonical surjective homomorphism
\begin{equation}\label{hom}
H^{l_0}_1(f^{-1}(t_0) ,\Z) \stackrel{\varphi}{\rightarrow} {\cal
M}_k( l_0,{\cal F}_\varepsilon )
\end{equation}
compatible with the action of $\pi_1(\C\setminus \Delta ,t_0)$.
The latter means that for every $\gamma \in \pi_1(\C\setminus
\Delta ,t_0)$ the diagram
\begin{eqnarray*}
  H^{l_0}_1(f^{-1}(t_0) ,\Z) & \stackrel{\varphi}{\rightarrow} & {\cal M}_k(
l_0,{\cal F}_\varepsilon ) \\
  \downarrow \gamma_* & &  \downarrow \gamma_* \\
H^{l_0}_1(f^{-1}(t_0) ,\Z) & \stackrel{\varphi}{\rightarrow} &
{\cal M}_k( l_0,{\cal F}_\varepsilon )
\end{eqnarray*}
commutes ($\gamma_*$ is the automorphism induced by $\gamma$).
Therefore $Ker (\varphi)$ is a subgroup of $H^{l_0}_1(f^{-1}(t_0)
,\Z)$, invariant under the action $\pi_1(\C\setminus \Delta
,t_0)$, and hence (\ref{repr}) is isomorphic to the induced
representation
$$
\pi_1(\C\setminus \Delta ,t_0) \rightarrow H^{l_0}_1(f^{-1}(t_0)
,\Z)/ Ker (\varphi)
$$
which is a subrepresentation of
$$
\pi_1(\C\setminus \Delta ,t_0) \rightarrow H^{l_0}_1(f^{-1}(t_0)
,\Z)^* .
$$
{\bf Proof of Theorem \ref{th1}.}
First of all, note that if $l_1,l_2 \in \pi_1(f^{-1}(t_0),P_0)$
and
$$
{\cal P}_{l_1,{\cal F}_\varepsilon }(t)= t + M_k(l_1,{\cal
F}_\varepsilon ,t)\varepsilon^k + {\rm O}(\varepsilon^{k+1}),\;\;
{\cal P}_{l_2,{\cal F}_\varepsilon }(t)= t + M_k(l_2,{\cal
F}_\varepsilon ,t)\varepsilon^k + {\rm O}(\varepsilon^{k+1})
$$
then
$$
{\cal P}_{l_1,{\cal F}_\varepsilon } \circ {\cal P}_{l_2,{\cal
F}_\varepsilon }(t) = {\cal P}_{l_2\circ l_1,{\cal F}_\varepsilon
}(t)= t+(M_k(l_1,{\cal F}_\varepsilon ,t)+M_k(l_2,{\cal
F}_\varepsilon ,t))\varepsilon^k + {\rm O}(\varepsilon^{k+1})
$$
(the proof repeats the arguments of Proposition \ref{sigma}). It
follows that

\begin{equation}
\label{mk}
M_k(l_1\circ l_2,{\cal F}_\varepsilon ,t)=M_k(l_2\circ l_1,{\cal
F}_\varepsilon ,t)= M_k(l_1,{\cal F}_\varepsilon ,t)+M_k(l_2,{\cal
F}_\varepsilon ,t).
\end{equation}

 The generating
function $M_k(t)$ is locally analytic and multivalued on
$\C\setminus \Delta $. For every determination
$\gamma_*M_k(l_0,{\cal F}_\varepsilon ,t)$ of $M_k(l_0,{\cal
F}_\varepsilon ,t)$ obtained after an analytic continuation along
a closed loop $\gamma \in \pi_1(\C\setminus \Delta ,t_0)$ it holds
\begin{equation}
\label{eq3}
\gamma_*M_k(l_0,{\cal F}_\varepsilon ,t))=
M_k(\gamma_* l_0,{\cal F}_\varepsilon ,t)
\end{equation}
where $l_0$ is (by abuse of notation) a free homotopy class of
closed loops on $f^{-1}(t_0)$. Indeed, let $l(t)\subset f^{-1}(t)$
be a continuous family of closed loops, $l(t_0)=l_0$. For each
$\tilde{t}_0$ we may define a holonomy map $ {\cal
P}_{l(\tilde{t}_0),{\cal F}_\varepsilon }(t)$ analytic in a
sufficiently small disc centered at $\tilde{t}_0$. It follows from
the definition of the holonomy map, that if $\tilde{t}_0, t_0$ are
fixed sufficiently close regular values of $f$, then $ {\cal
P}_{l(\tilde{t_0}),{\cal F}_\varepsilon }(t)$ and $ {\cal
P}_{l(t_0),{\cal F}_\varepsilon }(t)$ coincide in some open disc,
containing $\tilde{t}_0, t_0$. The same holds for the
corresponding generating functions. This shows that the analytic
continuation of $M_k(t)=M_k(l(t_0),{\cal F}_\varepsilon ,t)$ along
an interval connecting $t_0$ and $\tilde{t}_0$ is obtained by
taking a continuous deformation of the closed loop $l(t_0)$ along
this interval. Clearly this property of the generating function
holds true even without the assumption that $\tilde{t}_0, t_0$ are
close and for every path connecting $\tilde{t}_0, t_0$. This
proves the identity (\ref{eq3}).

Formula (\ref{eq3}) shows that
$$
k(l_0,{\cal F}_\varepsilon )= k(\gamma_*l_0,{\cal F}_\varepsilon),
\quad \forall \gamma \in \pi_1(\C\setminus \Delta ,t_0) .
$$
Let $l\subset f^{-1}(t_0)$ be a closed loop representing an
equivalence class in $H^{l_0}_1(f^{-1}(t_0) ,\Z)$. Then (\ref{mk})
implies that $k(l,{\cal F}_\varepsilon )\geq k(l_0,{\cal
F}_\varepsilon )$ and we define
$$
\varphi (l) = \left\{
\begin{array}{ccc}
 M_k(l,{\cal F}_\varepsilon ,t), & & \mbox{ \rm   if   }
k(l,{\cal
F}_\varepsilon ) = k(l_0,{\cal F}_\varepsilon ) \\
 0,& & \mbox{ \rm   if   } k(l,{\cal
F}_\varepsilon ) > k(l_0,{\cal F}_\varepsilon )\\
\end{array}
\right.
$$
\vspace{1ex}
Using the definitions of the abelian groups $H^{l_0}_1(f^{-1}(t_0)
,\Z) $ and $ {\cal M}_k( l_0,{\cal F}_\varepsilon )$ and the
identities (\ref{mk}), (\ref{eq3}), it is straightforward to check
that
\begin{itemize}
    \item $\varphi$ depends on the equivalence class of the loop
    $l$ in $H^{l_0}_1(f^{-1}(t_0)
,\Z) $;
    \item $\varphi(l)$ belongs to $ {\cal M}_k( l_0,{\cal F}_\varepsilon
    )$;
    \item $\varphi$ defines a surjective homomorphism (\ref{hom}) which
     is compatible with the action of $\pi _1(\C\setminus \Delta
,t_0)$ on $H^{l_0}_1(f^{-1}(t_0) ,\Z) $ and $ {\cal M}_k(
l_0,{\cal F}_\varepsilon )$.
\end{itemize}
Theorem \ref{th1} is proved. $\square$

\subsection{Main result}
\label{mainresult}
Our main result in this paper is the following.
\begin{theorem}
\label{main}
$ $
\begin{enumerate}
\item
If $H^{l_0}_1(f^{-1}(t_0) ,\Z) $ is of finite dimension, then the
generating function
 $M_k(t)=M_{k}(l_0,{\cal F}_\varepsilon, t)$ satisfies a
 linear differential equation
\begin{equation}
a_n(t)x^{(n)} + a_{n-1}(t) x^{(n-1)}+...+a_1(t) x' +a_0(t) x= 0
\label{fuchs}
\end{equation}
where $n \leq  dim\, H^{l_0}_1(f^{-1}(t_0) ,\Z)$ and $a_i(t)$ are
suitable analytic functions on $\C\setminus \Delta $.
\item
If, moreover, $M_{k}(t)$ is a function of moderate growth at any
$t_i \in \Delta$ and at $t=\infty$, then $(\ref{fuchs})$ is an
equation of Fuchs type.
\item
If in addition to the preceding hypotheses, the canonical map
\begin{equation}
H^{l_0}_1(f^{-1}(t_0) ,\Z) \rightarrow H_1(f^{-1}(t_0),\Z)
\label{injection}
\end{equation}
is injective, then $M_{k}(t)$ is an Abelian integral
\begin{equation}\label{abelian}
M_k(t)= \int_{l(t)} \omega
\end{equation}
where $\omega$ is a rational one-form on $\C^2$
and $l(t) \subset f^{-1}(t)$ is a continuous family of closed
loops, $l(t_0)=l_0$.

\end{enumerate}
\end{theorem}
{\bf Remarks.}
\begin{enumerate}
\item
Recall that a multivalued locally analytic function $g : \C
\setminus \Delta \rightarrow \C $ is said to be of moderate growth
if for every $\varphi_0 >0$ there exist constants $C, N>0$ such
that
$$
sup \{ |g(t)t^N|: 0< |t-t_i|< C,\;\; \mbox{\it Arg}\,(t-t_i) <
\varphi_0,\;\; t_i \in \Delta \} < \infty$$ and
$$
sup \{ |g(t)t^{-N}|: |t |> 1/C , \mbox{\it Arg}\,|t |  < \varphi_0
\} < \infty .
$$
\item

When (\ref{injection}) is not injective, the generating function could
still be an Abelian integral. Of course,
this depends on the unfolding ${\cal F}_\varepsilon $.
\item
If the dimension of $H^{l_0}_1(f^{-1}(t_0) ,\Z)$ is finite, we may
also suppose that (\ref{fuchs}) is irreducible. This makes
(\ref{fuchs}) unique (up to a multiplication by analytic
functions). The monodromy group of this equation is a subgroup of
the monodromy group associated to $l_0$, see  (\ref{monodromy1}).
It is clear that $M_{k}(t)$ may satisfy other equations with
non-analytic coefficients on $\C \setminus \Delta $.
\end{enumerate}
\vspace{2ex} \noindent {\bf Proof of Theorem \ref{main}.} Suppose
that $H^{l_0}_1(f^{-1}(t_0) ,\Z)$ is of finite dimension. Then
${\cal M}_k( l_0,{\cal F}_\varepsilon )=  H^{l_0}_1(f^{-1}(t_0)
,\Z)/ Ker(\varphi)$ is of finite dimension too, and let $g_i(t) =
M_k(l_i,{\cal F}_\varepsilon ,t) $, $i=1,...,n$ be a basis of the
{\it complex vector space} $V$ generated by ${\cal M}_k( l_0,{\cal
F}_\varepsilon )$, $\dim_{\C} V \leq \dim_{\Z} {\cal M}_k(
l_0,{\cal F}_\varepsilon )$. There is a  unique linear
differential equation of order $\dim_{\C} V$ satisfied by the
above generating functions (and hence by $M_k(l_0,{\cal
F}_\varepsilon ,t) $) having the form (\ref{fuchs}) which can be
equivalently written as
\begin{equation}
\label{fff} \det \pmatrix{ g_1 & g_1 '& ... &  g_1^{(n)} \cr g_2 &
g_2 '& ... &  g_2^{(n)} \cr
 . & . & ... &. \cr
 g_n & g_n '& ... &  g_n^{(n)} \cr \cr x & x'& ... & x^{(n)}}=0.
\end{equation}
The functions $g_1 ,g_2 ,...,  g_n $ are linearly independent over
$\C$ and define a complex  vector space invariant under the action
of $\pi _1(\C\setminus \Delta ,t_0)$. For a given $\gamma \in \pi
_1(\C\setminus \Delta ,t_0)$, let $\gamma _* \in Aut (V)$ be the
 automorphism  (\ref{monodromy1}) and denote (by abuse of notation) by
$\gamma_*a_i(t)$ the analytic
continuation of $a_i(t)$ along the loop $\gamma $.
The explicit form of the coefficients $a_i(t)$ as determinants (see
(\ref{fff})) implies that
 $\gamma_*a_i(t) = det(\gamma_* ) a_i(t)$. Therefore
$\gamma_*[ a_i(t)/a_n(t)]= a_i(t)/a_n(t)$, $a_i(t)/a_n(t)$ are
single-valued  and hence meromorphic functions on $\C\setminus
\Delta $. This proves the first claim of the theorem.  If in
addition $M_{k}(t)$ is of moderate growth, then  $g_i(t)$ are of
moderate growth too, $a_i(t)/a_n(t)$ are rational functions, and
the equation (\ref{fuchs}) is of Fuchs type (eventually with
apparent singularities).

Suppose finally that (\ref{injection}) is injective, which implies
that $H^{l_0}_1(f^{-1}(t_0) ,\Z)$ is a subgroup of the homology
group $ H_1(f^{-1}(t_0),\Z)$. By the algebraic de Rham theorem
\cite{gh} the first cohomology group of $f^{-1}(t_0)$ is generated
by polynomial one-forms. In particular, the dual space of
$H^{l_0}_1(f^{-1}(t_0) ,\Z)$ is generated by polynomial one-forms
$\omega _1,\omega _2,...,\omega _n$. Let $l_1(t),
l_2(t),...,l_n(t),l(t)\subset f^{-1}(t)$ be a continuous family of
closed loops, such that $l_1(t_0), l_2(t_0),...,l_n(t_0)$ defines
a basis of $ H^{l_0}_1(f^{-1}(t_0) ,\Z$, $l(t_0)=l_0$. The
determinant
\begin{equation}
\label{ddd} \det \pmatrix{ g_1 & \int_{l_1} \omega _1& \int_{l_1}
\omega _2&... & \int_{l_1}\omega _n \cr g_2 &\int_{l_2}\omega _1&
\int_{l_2} \omega _2&... & \int_{l_2}\omega _n \cr . & . & ... &.
\cr g_n &\int_{l_n}\omega _1& \int_{l_n} \omega _2&... &
\int_{l_n}\omega _n \cr M_{k} &\int_{l}\omega _1& \int_{l} \omega
_2&... & \int_{l}\omega _n }=0.
\end{equation}
developed with respect to the last row gives
$$
\alpha _0 M_k + \alpha _1 \int_{l}\omega _1 + \alpha_2\int_{l}
\omega _2 +...\alpha _n \int_{l}\omega _n = 0 .
$$
As $ H^{l_0}_1(f^{-1}(t_0) ,\Z)\subset H_1(f^{-1}(t_0),\Z)$ is
invariant under the action of $\pi _1(\C\setminus \Delta ,t_0)$,
then we deduce in the same way as before that $\alpha _i(t)/\alpha
_0(t)$ are rational functions. This completes the proof of the
theorem. $ \square$

We conclude the present section with some open questions.
Let $l_0(t)\subset f^{-1}(t)$ be a continuous family of ovals defined by
the real polynomial $f\in \R[x,y]$.

\vspace{2ex}
{\bf Open questions}
\begin{enumerate}
\item
Is it true that  the abelian group $H^{l_0}_1(f^{-1}(t_0) ,\Z)$ is
free, torsion free, finitely generated, or even stronger, $\dim
H^{l_0}_1(f^{-1}(t_0) ,\Z) \leq  dim H_1(f^{-1}(t_0),\Z)$? If not,
give counter-examples.
\item
Is it true that every generating function  of a polynomial
deformation ${\cal F}_\varepsilon $ of $df=0$ is of moderate
growth at any point $t\in \Delta $ or $t\in \infty$?
\item
Is it true that the monodromy representation (\ref{monodromy1})
has the following universal property : {\em for every $l \in
H^{l_0}_1(f^{-1}(t_0) ,\Z)$ there exists a polynomial deformation
${\cal F}_\varepsilon $ of $df=0$, such that the corresponding
generating function $\varphi(l)$ is not identically zero.}
 If this were true it would  imply that $H^{l_0}_1(f^{-1}(t_0) ,\Z)$ is
torsion-free, and whenever (\ref{injection}) is not injective,
then there exists a polynomial unfolding with corresponding
generating function which is not an Abelian integral of the form
(\ref{abelian}).
%Equivalently, this would mean that $G_{l_0}$ is the space dual to
%$H^{l_0}_1(f^{-1}(t_0) ,\Z)$ giving thus
%an analogue of the de Rham theorem.
\item Suppose that the canonical homomorphism
(\ref{canonical}) is surjective. Is it true that it is also
injective? Note that a negative answer would imply that the
representation (\ref{monodromy1}) is not universal (in the sense
of the preceding question). Indeed, if (\ref{canonical}) is
surjective, then  the orbit ${\cal O}_{l_0}$ generates the
homology group, and hence the generating function is always an
Abelian integral. The kernel of the canonical map
(\ref{canonical}) consist of free homotopy classes (modulo an
equivalence relation) homologous to zero, along which every
Abelian integral vanishes.
\end{enumerate}

\section{Examples}
\label{examples} In this section we show that the claims of
Theorem \ref{main} are non-empty. Namely, we apply it to
polynomial deformations $f$ of the simple singularities $y^2+x^4$,
$xy(x-y)$ of type $A_3$, $D_4$ respectively (see \cite[vol.
1]{avg} for this terminology). For a given loop $\delta(t) \subset
f^{-1}(t)\subset \C^2$ we shall compute the group
$H_1^\delta(f^{-1}(t) ,\Z)$. As  the abelian groups
$H_1^\delta(f^{-1}(t) ,\Z)$  are isomorphic, then when the choice
of $t$ is irrelevant we shall omit it. The same convention will be
applied to the cycles or closed loops on the fibers $f^{-1}(t)$.
 An equivalence class of loops in $H_1^\delta(f^{-1}(t_0)_t ,\Z) $ will be represented by a free
homotopy class of loops on $f^{-1}(t)$. Two such free homotopy
classes $\delta _1,\delta _2$ are composed in the following way:
take any two representative of $\delta _1,\delta _2$ in the
fundamental group of the surface $f^{-1}(t)$ and compose them.
This operation is compatible with the group law in
$H_1^\delta(f^{-1}(t_0)_t ,\Z) $, provided that $\delta _1,\delta
_2$ represent equivalence classes in it. The operation defines a
unique element in $H_1^\delta(f^{-1}(t_0)_t ,\Z) $ (represented
once again by a non-unique free homotopy class of loops).

\subsection{The $A_3$ singularity}
\begin{figure}
\begin{center}
\epsfig{file=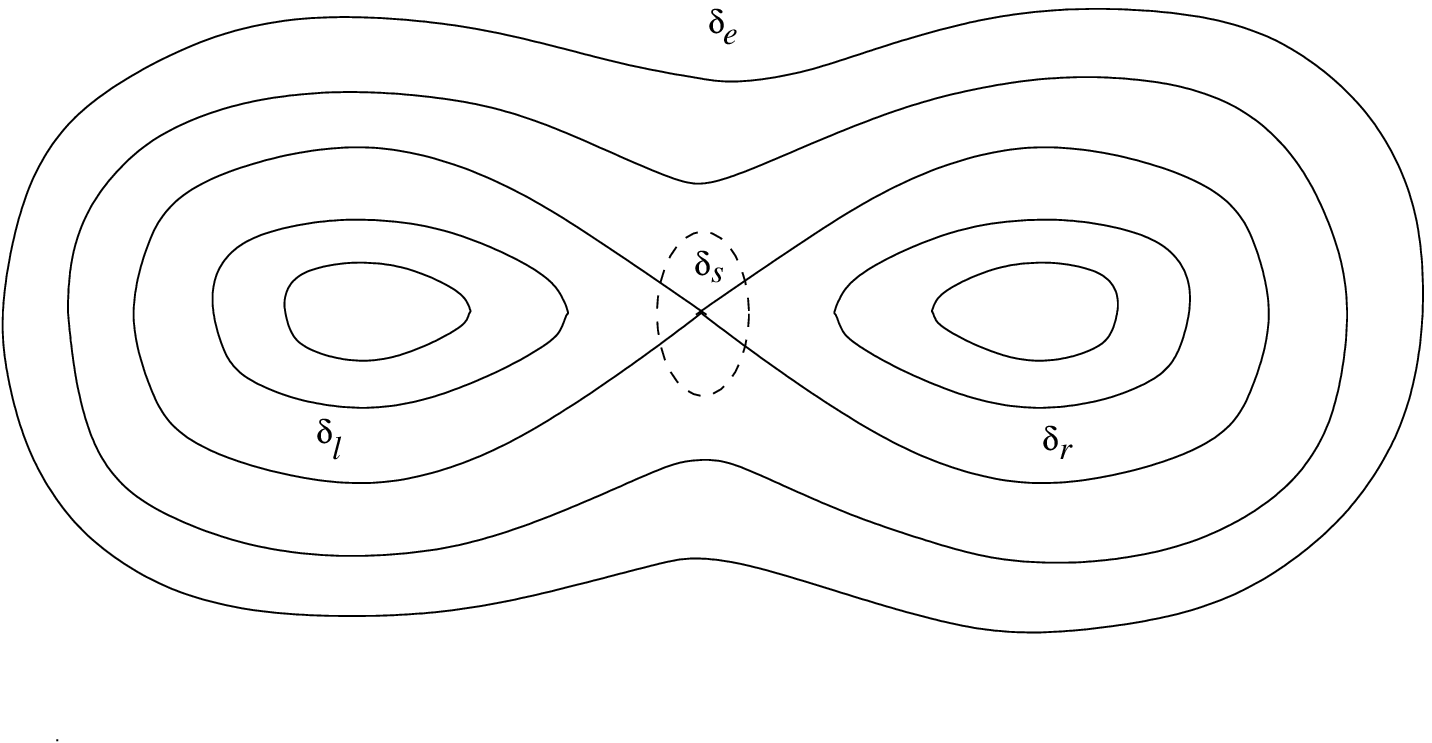}
\end{center}
\caption{The continuous families of ovals $ \delta _l, \delta _r$ and
$\delta _e$.}
\label{8}
\end{figure}

Take
$$
f(x,y)= \frac{y^2}{2} + \frac{(x^2-1)^2}{4}
$$
and denote by $\delta_e (t), \delta _l(t),\delta _r(t)$ respectively
the exterior,  left interior and right interior continuous family of
ovals defined by $\{(x,y)\in \R^2: f(x,y)=t\}$, see Fig. \ref{8}.
We denote by the same letters the corresponding continuous families of
free homotopy classes of loops defined on the universal covering space of
$\C \setminus \{0,1/4\}$, and fix $t_0 \neq 0$.
\begin{proposition}
\label{pr1}
We have
$$
H^{\delta _l}_1(f^{-1}(t_0),\Z)= H^{\delta
_r}_1(f^{-1}(t_0),\Z)=H_1(f^{-1}(t_0),\Z)=\Z^3,\quad H_{\delta_e
(t_0)}=\Z^2
$$
and the canonical map $H^{\delta _e}_1(f^{-1}(t_0),\Z) \rightarrow
H_1(f^{-1}(t_0),\Z)$ is injective.
\end{proposition}
Applying Theorem \ref{main} we get
\begin{corollary}
For every polynomial unfolding ${\cal F}_\varepsilon $ the generating function
$M_{\delta _e(t_0)}$ is an Abelian integral, provided that this function is
of moderate growth.
\end{corollary}
It is possible to show that $M_{\delta _e(t_0)}$ is always of moderate
growth (this will follow
from the explicit computations below). As for $M_{\delta _l(t_0)}$ and
$M_{\delta _r(t_0)}$, it
follows from \cite{gav99}  that these functions are always Abelian
integrals.

\begin{figure}
\begin{center}
\epsfig{file=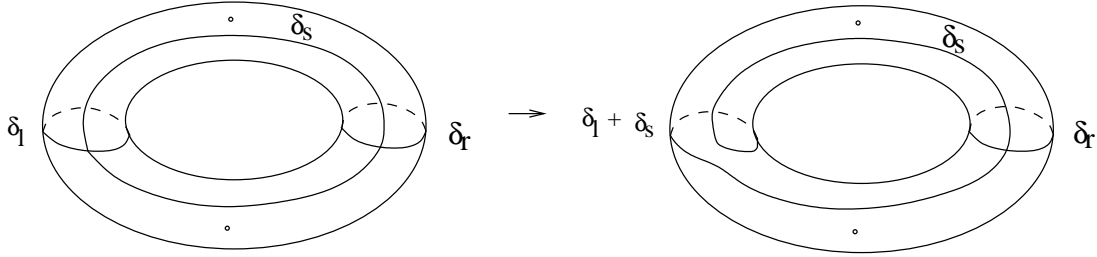}
\end{center}
\caption{The Dehn twist along the closed loop $\delta _l$.}
\label{dn}
\end{figure}

\begin{figure}
  \begin{center}
  \epsfig{file=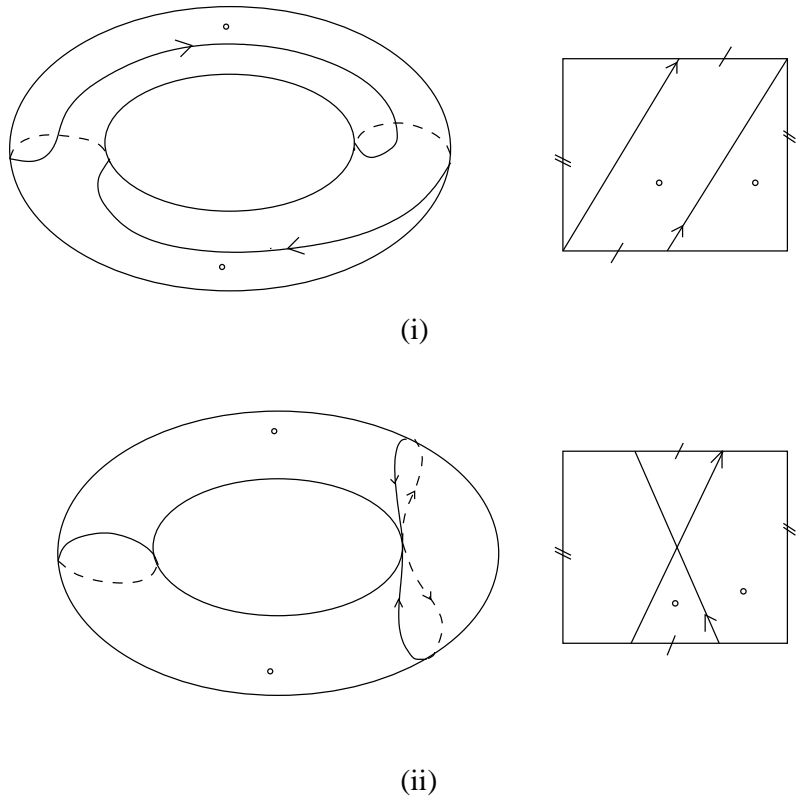}
  \end{center}
  \caption{ (i) The loop $l_{0*} \delta _s(t_0)$, and (ii)
  its first variation $\mbox{\it  Var}_{l_0}\delta _s(t_0)$.}
  \label{g2}
  \end{figure}

\vspace{2ex} \noindent {\bf Proof of Proposition \ref{pr1}.} The
affine curve $f^{-1}(t_0)$ is a torus with two removed points, and
hence $H_1(f^{-1}(t_0),\Z)=\Z^3$. We compute first $H^{\delta
_l}_1(f^{-1}(t_0),\Z)$. Let $t_0\in (0,1/4)$ and let $\delta
_s(t)\subset f^{-1}(t)$, $t\in (0,1/4)$, be the continuous family
of ``imaginary" closed loops (the ovals of $\{ y^2/2 - (x^2-1)^2/4
= t \}$) which tend to the saddle point $(0,0)$ as $t$ tends to
$1/4$. As before we denote by the same letter the  continuous
family of free homotopy classes of loops defined on the universal
covering space of $\C \setminus \{0,1/4\}$, and fix $t_0\neq
0,1/4$. Let $l_0,l_{1/4} \in \pi _1(\C\setminus \{0,1/4\},t_0)$ be
two simple loops making one turn about $0$ and $1/4$ respectively
in a positive direction. The group $\pi _1(\C\setminus
\{0,1/4\},t_0) $ acts on $\pi _1(f^{-1}(t_0))$ as follows. To the
loop $l_{1/4}$ corresponds an automorphism of $f^{-1}(t_0)$ which
is a Dehn twist along $\delta _s(t_0)$. Recall that a Dehn twist
of a surface along a closed loop is a diffeomorphism which is the
identity, except in a neighborhood of the loop. In a neighborhood
of the loop the diffeomorphism is shown on
 Fig. \ref{dn}, see \cite{poenaru}. The
usual Picard-Lefschetz formula \cite{avg} describes an automorphism
of the homology group induced by a Dehn twist along a ``vanishing" loop.
Therefore $l_{1/4*} \delta _s = \delta _s$ and $l_{1/4*} \delta
_l$ is the loop shown on Fig. \ref{dn}. We may also compose the
loops $\delta _s, l_{1/4*} \delta _l$ in the way explained in the
beginning of this section. The result is an equivalence class in
$H^{\delta _l}_1(f^{-1}(t_0),\Z)$ represented in a non-unique way
by a closed loop. The equivalence class $ \mbox{\it
Var}_{l_{1/4}}\delta _l= (l_{1/4}-id)_*\delta _l $ equals
therefore to the class represented by $\delta _s$, and hence
$\mbox{\it Var}_{l_{1/4}}^2\delta _l$ represents the zero class.
In a similar way we compute
 $l_{0*} \delta _s(t_0)$ which equals $\delta _s + \delta _r + \delta _l$,
as well as its first variation
$
\mbox{\it Var}_{l_0}\delta _s= (l_0-id)_*\delta _s
$
which equals $\delta _r + \delta _l$, see
 Fig. \ref{g2}. It follows that the second variation
 $\mbox{\it Var}^2_{l_0}\delta _s$ of $\delta _s$
may be represented by a loop homotopic to a point. We conclude
that $H^{\delta _l}_1(f^{-1}(t_0),\Z)$ is generated by equivalence
classes represented by $\delta _l, \delta _s, \delta _r$ and hence
it coincides with $H_1(f^{-1}(t_0),\Z)$ (generated by the same
loops). The computation of $H^{\delta _r}_1(f^{-1}(t_0),\Z)$ is
analogous.

\begin{figure}
\begin{center}
\epsfig{file=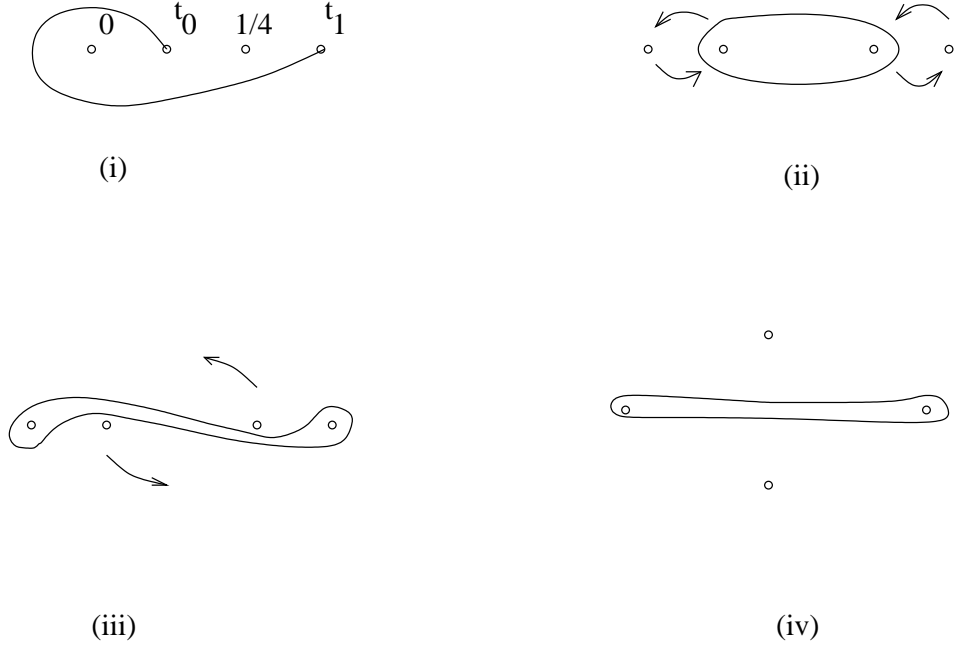}
\end{center}
\caption{The closed loops $\delta_e$ and $\delta_s$.}
\label{def}
\end{figure}
To compute $H^{\delta _e}_1(f^{-1}(t_0),\Z)$ we note that this
group coincides with $H^{\delta _s}_1(f^{-1}(t_0),\Z)$. Indeed,
take a loop $l\subset \C$ starting at $t_0\in (0,1/4)$ and
terminating at some $t_1\in (1/4,\infty)$ as it is shown on Fig.
\ref{def}. This defines a continuous family of (free homotopy
classes of) loops $\delta _s(t)$ along $l$. Then it follows from
Fig. \ref{def} that
$$
\delta _s(t_0)= \delta _e(t_1) \;
$$
and hence $H^{\delta _s(t_0)}_1(f^{-1}(t_0),\Z)= H^{\delta
_e(t_1)}_1(f^{-1}(t_0),\Z)$.  The loop $l_{0*} \delta _s(t_0)$ and
its first variation $ \mbox{\it Var}_{l_0}\delta _s(t_0)=
(l_0-id)_*\delta _s(t_0) $ were already computed (Fig. \ref{g2})
 and the second variation $\mbox{\it Var}^2_{l_0}\delta _s(t_0)$
may be represented by a loop homotopic to a point. Further,
$l_{1/4*} \delta _s(t_0) = \delta _s(t_0)$, and the first
variation $\mbox{\it Var}_{l_{1/4}}\mbox{\it Var}_{l_0}\delta
_s(t_0)$ of $\mbox{\it Var}_{l_0}\delta _s(t_0)$ along $l_{1/4}$
is a composition of free homotopy classes of $\delta _s$ (two
times). It follows that $H^{\delta _s}_1(f^{-1}(t_0),\Z)$ is
generated by $\delta _s$ and $\mbox{\it Var}_{l_0}\delta _s$. As
these loops are homologically independent we conclude that
$$H^{\delta _s}_1(f^{-1}(t_0),\Z) \rightarrow H_1(f^{-1}(t_0),\Z)$$
is injective and $H^{\delta _s}_1(f^{-1}(t_0),\Z)=\Z^2$. The
proposition is proved. $\square$

\subsection{Calculation of the generating function in the $A_3$ case}
In what follows we compare the above geometric approach to the
combinatorial approach based on Fran\c{c}oise's      recursion formulae.
We shall prove a stronger result allowing us to set up an explicit
upper bound to the number of zeros in $\Sigma$ of the displacement map
${\cal P}_\varepsilon (t)-t$ for small $\varepsilon$.
Below we use the standard notation $H$ of the Hamiltonian function,
$$H=\frac{y^2}{2}+\frac{(x^2-1)^2}{4}.$$

We say that $A$ is a polynomial of weighted degree $m$ in $x,y,H$
provided that
$$A(x,y,H)=\sum_{i+j+2k\leq m}a_{ijk}x^iy^jH^k$$
(namely, the weight of $x,y$ is one and the weight of $H$ is assumed to
be two). Clearly, a polynomial in $x,y$ allows a representation through
different weighted polynomials in $x,y,H,$ possibly of different weighted
degrees, depending on the way the powers $x^i$ with $i>3$ were expressed.
However, any polynomial has a unique representation through a weighted
polynomial in a normal form which means that the latter contains powers
$x^i$ with $i\leq 3$ only. We will not assume that the weighted polynomials
we consider bellow are taken in a normal form.

Set $\sigma_k=x^kydx$ and $I_k(t)=\int_{\delta(t)}\sigma_k$, $k=0,1,2,$
where $\delta(t)$ is an oval contained in the level set $\{H=t\}$.

%\vspace{2ex}
%\noindent
\begin{proposition}
\label{petrov}
For any one-form $\omega_m=A_m(x,y,H)dx+B_m(x,y,H)dy$ with polynomial
coefficients of weighted degree $m$, the following decomposition holds:
\begin{equation}
\label{1}
\omega_m=dG_{m+1}(x,y,H)+g_{m-1}(x,y,H)dH+\alpha_{m-1}(H)\sigma_0
+\beta_{m-2}(H)\sigma_1 +\gamma_{m-3}(H)\sigma_2
\end{equation}
where $G_k,g_k,\alpha_k,\beta_k,\gamma_k$ are polynomials in their
arguments of weighted degree $k$.
\end{proposition}
%
%\vspace{2ex}
%\noindent
Below, we will denote by $\alpha_k,\beta_k,\gamma_k$
polynomials of weighted degree $k$ in $H$, by $G_k, g_k$
polynomials of weighted degree $k$ in $x,y,H$,
and by $\omega_k$ one-forms with polynomial coefficients of weighted
degree $k$ in $x,y,H$. (Possibly, different polynomials and one-forms
of the same degree and type will be denoted by the same letter.)

\vspace{2ex}
\noindent
{\bf Proof of Proposition \ref{petrov}.}
The proof is similar to the proof of Lemma 1
in \cite{ili00}
which concerned the elliptic case $H=\frac12y^2+\frac12x^2-\frac13x^3$.
It is sufficient to consider the case when the coefficients of the one-form
do not depend
on $H$. As in \cite{ili00}, one can easily see that the problem reduces to
expressing the one-forms $y^jdx, xy^jdx, x^2y^jdx$ in the form (\ref{1}). We
have
$$\begin{array}{rl}
y^jdx&{\displaystyle \!\!\!=\frac{4j}{2j+1}Hy^{j-2}dx
+\frac{j}{2j+1}(x^2-1)y^{j-2}dx
-\frac{j}{2j+1}xy^{j-2}dH+d\frac{xy^j}{2j+1},}\\[3mm]
xy^jdx&{\displaystyle \!\!\!=\frac{2j}{j+1}Hxy^{j-2}dx
-\frac{j}{2j+2}(x^2-1)y^{j-2}dH+ d\frac{(x^2-1)y^j}{2j+2},}\\[3mm]
x^2y^jdx&{\displaystyle \!\!\!-\,\frac{1}{2j+3}y^jdx=
\frac{4j}{2j+3}Hx^2y^{j-2}dx
-\frac{j}{2j+3}(x^3-x)y^{j-2}dH+d\frac{(x^3-x)y^j}{2j+3}.}
\end{array}$$
>From the second equation we obtain immediately that
$xy^jdx=c_jH^{\frac{j-1}{2}}\sigma_1+dG_{j+2}+g_jdH$
($c_j=0$ for $j$ even, $c_j>0$ for $j$ odd) which yields
$$xA_{m-1}(y)dx=\beta_{m-2}(H)\sigma_1+dG_{m+1}+g_{m-1}dH.$$
Taking notation $\theta_j=(y^jdx, x^2y^jdx)^\top$,
$\Theta_j=(dG_{j+1}+g_{j-1}dH, dG_{j+3}+g_{j+1}dH)^\top$, one can
rewrite the system formed by the first and the third equation above
in the form
$$\theta_j=\Lambda_j(H)\theta_{j-2}+\Theta_j, \quad
\Lambda_j(H)=\frac{j}{2j+1}\left(\begin{array}{cc}
4H-1 &1\\[2mm]
{\displaystyle\frac{4H-1}{2j+3}}&{\displaystyle\frac{4(2j+1)H+1}{2j+3}}
\end{array}\right).$$
As $\Lambda_j\Theta_{j-2}=\Theta_j$, this implies that
$\theta_j=\Lambda_j\Lambda_{j-2}\ldots\Lambda_3\theta_1 +\Theta_j$ for $j$ odd
and $\theta_j=\Lambda_j\Lambda_{j-2}\ldots\Lambda_2\theta_0 +\Theta_j$ for $j$
even, which in both cases is equivalent to
\begin{equation}
\label{suppl}
\begin{array}{l}
y^jdx=\alpha_{j-1}(H)\sigma_0+\gamma_{j-3}(H)\sigma_2+ dG_{j+1}+g_{j-1}dH,
\\[2mm]
x^2y^jdx=\alpha_{j-1}(H)\sigma_0+\gamma_{j-1}(H)\sigma_2+ dG_{j+3}+g_{j+1}dH
\end{array}
\end{equation}
where the coefficients at $\sigma_0, \sigma_2$ vanish for $j$ even.
Applying the last two relations with $j\leq m$ and $j\leq m-2$ respectively,
we obtain the result. $\Box$

\vspace{2ex}
\noindent
The above decomposition (\ref{1}) is the basic tool for calculating the
generating functions. For the two period annuli inside the eight-loop
(level sets $t\in (0,\frac14)$), one has
$$\int_{\delta(t)}\omega_m=\alpha_{m-1}(t)I_0(t)+\beta_{m-2}(t)I_1(t)
+\gamma_{m-3}(t)I_2(t),$$
and for $0<t<\frac14$,
$$\int_{\delta(t)}\omega_m\equiv 0\;\;\Leftrightarrow\;\;
\alpha_{m-1}(t)=\beta_{m-2}(t)=\gamma_{m-3}(t)\equiv 0\;\;\Leftrightarrow\;\;
\omega_m=dG_{m+1}+g_{m-1}dH.$$
This means that the internal period annuli satisfy the so called $(*)$ property
\cite{fra96} and the generating functions are determined from the integration
of polynomial one-forms calculated in a recursive procedure. More explicitly,
consider a small polynomial perturbation
\begin{equation}
\label{system}
\begin{array}{l}
\dot{x}=H_y+\varepsilon f(x,y),\\
\dot{y}=-H_x+\varepsilon g(x,y),\end{array}
\end{equation}
which can be rewritten as $dH-\varepsilon\omega_n=0$
with $\omega_n=g(x,y)dx-f(x,y)dy$ and $n$ the degree of the perturbation.
Then in $(0,\frac14)$, the first nonzero generating function is given by
$$M_k(t)=\int_{\delta(t)}\Omega_k,\;\;\mbox{\rm where}\;\;
\Omega_1=\omega_n,\;\; \Omega_k=q_{k-1}\Omega_1\;\;\mbox{\rm and}\;\;
\Omega_{k-1}=dQ_{k-1}+q_{k-1}dH.$$
Making use of (\ref{1}), it is then easily seen by induction that
$q_{k-1}$ is a polynomial of weighted degree $(k-1)(n-1)$, therefore
$\Omega_k$ is a polynomial one-form of weighted degree
$m=k(n-1)+1$ which proves that
\begin{equation}
\label{kth}
M_k(t)=\alpha_{[\frac{k(n-1)}{2}]}(t)I_0(t)
+\beta_{[\frac{k(n-1)-1}{2}]}(t)I_1(t)
+\gamma_{[\frac{k(n-1)-2}{2}]}(t)I_2(t)
\end{equation}
where $\alpha_j,\beta_j,\gamma_j$ are polynomials in $t$ of degree at most $j$.

For the period annulus outside the eight-loop
(level sets $t\in (\frac14,\infty)$), one has
$$\int_{\delta(t)}\omega_m=\alpha_{m-1}(t)I_0(t)+\gamma_{m-3}(t)I_2(t),$$
and for $\frac14<t<\infty$,
$$\int_{\delta(t)}\omega_m\equiv 0\;\;\Leftrightarrow\;\;
\alpha_{m-1}(t)=\gamma_{m-3}(t)\equiv 0 \;\;\Leftrightarrow\;\;
\omega_m=dG_{m+1}+g_{m-1}dH+\beta_{m-2}(H)\sigma_1$$
since $I_1(t)\equiv 0$ which is caused by symmetry of the oval.
Therefore the outer period annulus does not satisfy the $(*)$ property
which makes this case troublesome and we shall deal with it
until the end of this section.

Take a point $(x,y)$ lying on a certain level set $H=t$ for
a fixed $t>\frac14$ and let $(a,0)$ be the intersection point of the
level curve with the negative $x$-axis. Denote by $\delta(x,y)\subset\{H=t\}$
the oriented curve in the $(\xi,\eta)$ plane connecting $(a,0)$ and $(x,y)$
in a clockwise direction. Consider the function $\varphi$ determined
by the formula (see formula (2.5) in \cite{ili98})
$$\varphi(x,y)=\int_{\delta(x,y)}\frac{\xi d\xi}{\eta}.$$
As $I_1(t)=\int_{\delta(t)}xydx\equiv 0$, this is also true for
$I_1'(t)=\int_{\delta(t)}\frac{xdx}{y}$ which implies that
$\varphi(\pm a,0)=0$.
Therefore, $\varphi(x,y)$ is single-valued and hence an analytic function
in the domain outside the eight-loop.
In \cite{iliper}, $\varphi$ was expressed as
$$\varphi(x,y)=\frac{1}{\sqrt2}\left(
{\rm arctan}\frac{x^2-1}{y\sqrt2}-\frac{\pi}{2}{\rm sign}\, y\right)
=\frac{{\rm sign}\, y}{\sqrt2}\left(
{\rm arcsin}\frac{x^2-1}{2\sqrt{H}}-\frac{\pi}{2}\right).$$
In \cite{jmp02}, the authors expressed $\varphi$ by a complex logarithmic
function
$$\varphi=\frac{i}{2\sqrt2}\log\frac{x^2-1+i\sqrt2y}{x^2-1-i\sqrt2y}$$
and used in their proofs the properties of $\varphi$ on the corresponding
Riemann surface. The concrete expression of the function $\varphi$ is
inessential in our analysis. We will only make use of the identities
(\ref{2}) below and the fact that $\varphi$ there is determined up to an
additive constant, whilst the first nonvanishing generating function $M_k$
is independent on such a constant.

\vspace{2ex}
\noindent
Let us denote for short $G=\frac14(x^2-1)y.$
Using direct calculations, one can establish easily the following identities:
\begin{equation}
\label{2}
\begin{array}{l}
\sigma_1=xydx=dG+Hd\varphi,\\[2mm]
{\displaystyle Hd\varphi=\frac{xy}{2}dx-\frac{x^2-1}{4}dy},\\[3mm]
(x^2-1)d\varphi={\displaystyle\frac{y}{2H}}dH-dy,\\[3mm]
yd\varphi={\displaystyle xdx-\frac{x^2-1}{4H}dH},\\[3mm]
xd\varphi={\displaystyle \frac{(5x^2-1)y}{4H}dx-d\left(\frac{xG}{H}\right)
-\frac{xG}{H^2}dH.}
\end{array}
\end{equation}

\vspace{2ex}
\noindent
Making use of the first identity in (\ref{2}), we can rewrite (\ref{1}) as
\begin{equation}
\label{3}
\begin{array}{rl}
\omega_m&\!\!=dG_{m+1}+g_{m-1}dH+\beta_md\varphi+\alpha_{m-1}(H)\sigma_0
+\gamma_{m-3}(H)\sigma_2\\[2mm]
&\!\!=d(G_{m+1}+\varphi \beta_m)+(g_{m-1}-\varphi \beta_m')dH
+\alpha_{m-1}(H)\sigma_0 +\gamma_{m-3}(H)\sigma_2,
\end{array}
\end{equation}
with some new $G_k, g_k$ and $\beta_m$ satisfying $\beta_m(0)=0$.
%
%\vspace{2ex}
%\noindent
\begin{lemma}
\label{basic}
 For any nonnegative integer $l$ and one-form of weighted
degree $m\geq 0$, the following identity holds:
$$\varphi^l\omega_m=d\sum_{j=0}^{l+1}\frac{\varphi^j}{H^{l-j}}G_{m+3l-3j+1}
+\sum_{j=0}^{l+1}\frac{\varphi^j}{H^{l-j+1}}g_{m+3l-3j+1}dH$$
\vspace{-2ex}
$$ +\sum_{j=0}^l\frac{\varphi^j}{H^{l-j}}\alpha_{m+3l-3j-1}\sigma_0
+\sum_{j=0}^l\frac{\varphi^j}{H^{l-j}}\gamma_{m+3l-3j-3}\sigma_2.$$
\end{lemma}
%
%\vspace{2ex}
%\noindent
{\bf Proof.} By the first equation in (\ref{3}), we have
\begin{equation}
\label{3.5}
\begin{array}{rl}
\varphi^l\omega_m&{\displaystyle\!\!=d\left(\varphi^lG_{m+1}
+\frac{\varphi^{l+1}}{l+1}\beta_m\right)
+\left(\varphi^lg_{m-1}-\frac{\varphi^{l+1}}{l+1}\beta_m'\right)dH}\\[2mm]
&\!\!+\varphi^l\alpha_{m-1}\sigma_0
+\varphi^l\gamma_{m-3}\sigma_2-l\varphi^{l-1}G_{m+1}d\varphi.
\end{array}
\end{equation}
Using the second equation in (\ref{2}), we can rewrite this identity as
$$\varphi^l\omega_m=d(\varphi^lG_{m+1}+H\varphi^{l+1}G_{m-2})
+(\varphi^lg_{m-1}+\varphi^{l+1}g_{m-2})dH$$
\vspace{-2ex}
$$+\varphi^l\alpha_{m-1}\sigma_0
+\varphi^l\gamma_{m-3}\sigma_2+\frac{l}{H}\varphi^{l-1}\omega_{m+3}.$$
By iteration procedure, we get
$$\varphi^l\omega_m=\sum_{j=0}^l\frac{j!}{H^j}\pmatrix{l\cr j}
\left[d(\varphi^{l-j}G_{m+3j+1}+H\varphi^{l-j+1}G_{m+3j-2}) \right.$$
$$\left.+(\varphi^{l-j}g_{m+3j-1}+\varphi^{l-j+1}g_{m+3j-2})dH
+\varphi^{l-j}\alpha_{m+3j-1}\sigma_0
+\varphi^{l-j}\gamma_{m+3j-3}\sigma_2\right]$$
$$=\sum_{j=0}^l
\left[d\left(\frac{\varphi^{l-j}}{H^j}G_{m+3j+1}+\frac{\varphi^{l-j+1}}
{H^{j-1}}G_{m+3j-2}\right)
+\left(\frac{\varphi^{l-j}}{H^{j+1}}g_{m+3j+1}+
\frac{\varphi^{l-j+1}}{H^j}g_{m+3j-2}\right)dH\right.$$
$$\left. +\frac{\varphi^{l-j}}{H^j}\alpha_{m+3j-1}\sigma_0
+\frac{\varphi^{l-j}}{H^j}\gamma_{m+3j-3}\sigma_2\right]$$
$$=d\sum_{j=0}^{l+1}\frac{\varphi^j}{H^{l-j}}G_{m+3l-3j+1}
+\sum_{j=0}^{l+1}\frac{\varphi^j}{H^{l-j+1}}g_{m+3l-3j+1}dH$$
$$ +\sum_{j=0}^l\frac{\varphi^j}{H^{l-j}}\alpha_{m+3l-3j-1}\sigma_0
+\sum_{j=0}^l\frac{\varphi^j}{H^{l-j}}\gamma_{m+3l-3j-3}\sigma_2.\;\Box$$

\vspace{1ex}
\noindent
Unfortunately, one cannot use directly Lemma \ref{basic} to prove
Proposition \ref{two} and Theorem \ref{one}. Indeed, by the second
equation in (\ref{3}), we see that the function
$q_1$ is a first degree polynomial with respect to $\varphi$ which agrees
with Proposition \ref{two} for $k=1$. By applying Lemma \ref{basic},
we then conclude that $q_2$ would contain terms with denominators $H^2$,
which does not agree with Proposition \ref{two} when $k=2$.
The core of the problem is the following. Let us
express $\Omega_k$, the differential one-form used to calculate $M_k(t)$,
in the form $\Omega_k=dQ_k+q_kdH+a_k(H)\sigma_0+b_k(H)\sigma_2$.
Then $M_k(t)\equiv 0$ is equivalent to $a_k=b_k\equiv 0$.
However, the vanishing of $a_k$ and $b_k$ implies the vanishing of some
"bad" terms in $q_k$ as well. Without removing these superfluous terms in
$q_k$, one cannot derive the precise formulas of $M_{k+1}$ and $q_{k+1}$
during the next step. Hence, the precise result we are going to establish
requires much more efforts. The proof of our theorem therefore consists of
a multi-step reduction allowing us to detect and control these "bad" terms.
As the first step, we derive below some preliminary formulas.

Consider the function $G_{m+1}$ in formula (\ref{3.5}).
As it is determined up to an additive constant, one can write
$$G_{m+1}(x,y,H)=ax+(x^2-1)G_{m-1}(x)+yG_m(x,y)+HG_{m-1}(x,y,H)$$
which together with (\ref{2}) yields
$$-lG_{m+1}d\varphi=\omega_{m+1}+\frac{g_{m+2}}{H}dH+a_lxd\varphi,
\quad a_l=const. $$
Therefore, by (\ref{3.5}),
$$\varphi^l\omega_m=\varphi^{l-1}\omega_{m+1}
+d\left(\varphi^lG_{m+1}+\frac{\varphi^{l+1}}{l+1}\beta_m\right)
+\left(\frac{\varphi^{l-1}}{H}g_{m+2}+\varphi^lg_{m-1}
-\frac{\varphi^{l+1}}{l+1}\beta_m' \right)dH$$
$$+a_l\varphi^{l-1}xd\varphi +\varphi^l\alpha_{m-1}\sigma_0
+\varphi^l\gamma_{m-3}\sigma_2.$$
By iteration, one obtains
$$\varphi^l\omega_m= d\sum_{j=0}^l\left(\varphi^jG_{m+l-j+1}
+\frac{\varphi^{j+1}}{j+1}\beta_{m+l-j}\right)
+\sum_{j=0}^l\left(\varphi^jg_{m+l-j-1}
-\frac{\varphi^{j+1}}{j+1}\beta_{m+l-j}'\right)dH$$
\vspace{-2mm}
$$+\sum_{j=1}^{l}\frac{\varphi^{j-1}}{H}g_{m+l-j+2}dH
+\sum_{j=1}^la_j\varphi^{j-1}xd\varphi
+\sum_{j=0}^l\varphi^j\alpha_{m+l-j-1}\sigma_0
+\sum_{j=0}^l\varphi^j\gamma_{m+l-j-3}\sigma_2.$$
After a rearrangement, we get
\begin{equation}
\label{4}
\begin{array}{rl}
\varphi^l\omega_m&\!\!{\displaystyle = d\sum_{j=0}^{l+1}\varphi^jG_{m+l-j+1}
+\left(\sum_{j=0}^{l-1}\frac{\varphi^{j}}{H}g_{m+l-j+1}+\varphi^lg_{m-1}
-\varphi^{l+1}G_m'\right)dH}\\[5mm]
&\!\!{\displaystyle +\sum_{j=0}^{l-1}a_j\varphi^{j}xd\varphi
+\sum_{j=0}^l\varphi^j\alpha_{m+l-j-1}\sigma_0
+\sum_{j=0}^l\varphi^j\gamma_{m+l-j-3}\sigma_2}
\end{array}
\end{equation}
where $G_m=G_m(H)$ and $G_m(0)=0$. Using (\ref{4}), we then obtain
\begin{equation}
\label{5}
\begin{array}{rl}
{\displaystyle \frac{\varphi^l}{H}\omega_m} & \!\!{\displaystyle =
d\sum_{j=0}^{l+1}\frac{\varphi^j}{H}G_{m+l-j+1}
+\left(\sum_{j=0}^{l}\frac{\varphi^{j}}{H^2}g_{m+l-j+1}
-\varphi^{l+1}(G_m/H)'\right)dH}\\[5mm]
&\!\!{\displaystyle +\sum_{j=0}^{l-1}a_j\frac{\varphi^{j}}{H}xd\varphi
+\sum_{j=0}^l\frac{\varphi^j}{H}\alpha_{m+l-j-1}\sigma_0
+\sum_{j=0}^l\frac{\varphi^j}{H}\gamma_{m+l-j-3}\sigma_2.}
\end{array}
\end{equation}
More generally, for any $k\geq 2$,
\begin{equation}
\label{6}
\begin{array}{rl}
{\displaystyle \frac{\varphi^l}{H^k}\omega_m} & \!\!{\displaystyle =
d\sum_{j=0}^{l+1}\frac{\varphi^j}{H^k}G_{m+l-j+1}
+\left(\sum_{j=0}^{l}\frac{\varphi^{j}}{H^{k+1}}g_{m+l-j+1}
+\frac{\varphi^{l+1}}{H^k}\beta_{m-2}\right)dH}\\[5mm]
&\!\!{\displaystyle +\sum_{j=0}^{l-1}a_j\frac{\varphi^{j}}{H^k}xd\varphi
+\sum_{j=0}^l\frac{\varphi^j}{H^k}\alpha_{m+l-j-1}\sigma_0
+\sum_{j=0}^l\frac{\varphi^j}{H^k}\gamma_{m+l-j-3}\sigma_2.}
\end{array}
\end{equation}
After making the above preparation, take again a perturbation (\ref{system})
or equivalently $dH-\varepsilon\omega_n=0$ where $\omega_n$ is a polynomial
one-form in $(x,y)$ of degree $n$ and consider the related displacement map
(\ref{displ}).
%
%\vspace{2ex}
%\noindent
\begin{proposition}
\label{two}
 Assume that $M_1(t)=\ldots=M_k(t)\equiv 0$.
Then $\Omega_k=dQ_k+q_kdH$, with
\begin{equation}
\label{8a}
q_k=\sum_{j=0}^{k-1}\frac{\varphi^j}{H^{k-j-1}}g_{kn+k-3j-2}
+\varphi^kg_{k(n-2)}.
\end{equation}
\end{proposition}
%
%\vspace{2ex}
%\noindent
{\bf Proof.} The proof is by induction.  Assume that $q_k$ takes the form
(\ref{8a}), then $\Omega_{k+1}=q_k\Omega_1=q_k\omega_n$ can be written as
\begin{equation}
\label{ind}
\Omega_{k+1}=\sum_{l=0}^{k-1}\frac{\varphi^l}{H^{k-l-1}}
\omega_{(k+1)(n+1)-3l-3}+\varphi^k\omega_{(k+1)(n-2)+2}.
\end{equation}
Using (\ref{6}), we obtain that $M_{k+1}(t)=\int_{\delta(h)}\Omega_{k+1}^*$
where
$$\begin{array}{l}{\displaystyle
\Omega_{k+1}^*=\sum_{l=0}^{k-1}\left(\sum_{j=0}^{l-1}
a_{jl}\frac{\varphi^j}{H^{k-l-1}}xd\varphi +\sum_{j=0}^l
\frac{\varphi^j}{H^{k-l-1}}\alpha_{(k+1)(n+1)-j-2l-4}\sigma_0
\right.}\\[3mm]
{\displaystyle\left.\hspace{2ex}
+\sum_{j=0}^l\frac{\varphi^j}{H^{k-l-1}}\gamma_{(k+1)(n+1)-j-2l-6}\sigma_2
\right)}\\[3mm]
{\displaystyle\hspace{2ex}+\sum_{j=0}^{k-1} a_{jk}\varphi^jxd\varphi
+\sum_{j=0}^k \varphi^j\alpha_{(k+1)(n-2)+k-j+1}\sigma_0
+\sum_{j=0}^k\varphi^j\gamma_{(k+1)(n-2)+k-j-1}\sigma_2}\\[3mm]
{\displaystyle=\Omega_{k+1}^{**}+
\sum_{j=0}^{k-1} \varphi^j\frac{\alpha_{(k+1)(n+1)-3j-4}}{H^{k-j-1}}\sigma_0
+\sum_{j=0}^{k-1}\varphi^j\frac{\gamma_{(k+1)(n+1)-3j-6}}{H^{k-j-1}}
\sigma_2}\\[2mm]
{\displaystyle\hspace{2ex} + \varphi^k\alpha_{(k+1)(n-2)+1}\sigma_0
+ \varphi^k\gamma_{(k+1)(n-2)-1}\sigma_2}
\end{array}$$
and
$$\Omega_{k+1}^{**}=\sum_{l=0}^{k-2}\frac{\delta_{2k-2l-4}(H)}{H^{k-l-2}}
\varphi^lxd\varphi+a_{k-1}\varphi^{k-1}xd\varphi.$$
We now apply Lemma \ref{basic} (with $m=3$) to $\Omega_{k+1}^{**}$. Thus,
$$\Omega_{k+1}^{**}=\sum_{l=0}^{k-2}\frac{\delta_{2k-2l-4}}{H^{k-l-1}}
\varphi^l\omega_3+a_{k-1}\varphi^{k-1}xd\varphi$$
$$=\sum_{l=0}^{k-2}\frac{\delta_{2k-2l-4}}{H^{k-l-1}}
\sum_{j=0}^{l+1}\left[d\frac{\varphi^j}{H^{l-j}}G_{3l-3j+4}
+\frac{\varphi^j}{H^{l-j+1}}g_{3l-3j+4}dH\right]$$
$$+\sum_{l=0}^{k-2}\frac{\delta_{2k-2l-4}}{H^{k-l-1}}
\sum_{j=0}^{l}\frac{\varphi^j}{H^{l-j}}(\alpha_{3l-3j+2}\sigma_0
+\gamma_{3l-3j}\sigma_2)+a_{k-1}\varphi^{k-1}xd\varphi$$
$$=\sum_{l=0}^{k-2}\sum_{j=0}^{l+1}\left[d\frac{\varphi^j}
{H^{k-j-1}}G_{2k-3j+l}
+\frac{\varphi^j}{H^{k-j}}g_{2k-3j+l}dH\right]$$
$$+\sum_{l=0}^{k-2}\sum_{j=0}^{l}\frac{\varphi^j}{H^{k-j-1}}
(\alpha_{2k-3j+l-2}\sigma_0
+\gamma_{2k-3j+l-4}\sigma_2)+a_{k-1}\varphi^{k-1}xd\varphi$$
$$=\sum_{j=0}^{k-1}\left[d\frac{\varphi^jG_{3k-3j-2}}{H^{k-j-1}}
+\frac{\varphi^jg_{3k-3j-2}}{H^{k-j}}dH\right]$$
$$+\sum_{j=0}^{k-2}\frac{\varphi^j}{H^{k-j-1}}(\alpha_{3k-3j-4}\sigma_0
+\gamma_{3k-3j-6}\sigma_2)+a_{k-1}\varphi^{k-1}xd\varphi.$$
We have proved that
\begin{equation}
\label{9}
\begin{array}{rl}
\Omega_{k+1}^{*}=&\!\!
{\displaystyle  \sum_{j=0}^{k-1}\left[d\frac{\varphi^jG_{3k-3j-2}}{H^{k-j-1}}
+\frac{\varphi^jg_{3k-3j-2}}{H^{k-j}}dH\right]}\\[3mm]
&\!\!{\displaystyle+\sum_{j=0}^{k-1}\varphi^j\frac{\alpha_{(k+1)(n+1)-3j-4}}
{H^{k-j-1}}\sigma_0
+\sum_{j=0}^{k-1}\varphi^j\frac{\gamma_{(k+1)(n+1)-3j-6}}{H^{k-j-1}}
\sigma_2}\\[3mm]
&\!\!{\displaystyle + \varphi^k\alpha_{(k+1)(n-2)+1}\sigma_0
+ \varphi^k\gamma_{(k+1)(n-2)-1}\sigma_2+a_{k-1}\varphi^{k-1}xd\varphi.}
\end{array}
\end{equation}
We finish this step of the proof of Proposition \ref{two} by noticing
that if $M_{k+1}(t)=\int_{\delta(t)}\Omega_{k+1}^*\equiv 0$, then
the constant $a_{k-1}$ and the coefficients of all the polynomials
$\alpha_j,\gamma_j$ in (\ref{9}) are zero. The proof of this claim
is the same as the proof of Proposition \ref{three} below
and for this reason we omit it here. Therefore, equation
(\ref{9}) reduces to $\Omega^*_{k+1}=dQ_{k+1}^*+q_{k+1}^*dH$.

Next, applying to (\ref{ind}) the more precise identities (\ref{4}),
(\ref{5}) along with (\ref{6}), we see that
$\Omega_{k+1}=dQ_{k+1}+q_{k+1}dH+\Omega_{k+1}^*$ and
moreover, the coefficient at $dH$ is
$$q_{k+1}=\sum_{l=0}^{k-3}\left(\sum_{j=0}^l\frac{\varphi^j}{H^{k-l}}
g_{(k+1)(n+1)-j-2l-2}+\frac{\varphi^{l+1}}{H^{k-l-1}}g_{(k+1)(n+1)-3l-5}
\right)$$
$$+\sum_{j=0}^{k-2}\frac{\varphi^j}{H^2}
g_{(k+1)(n-2)+k-j+5}+\varphi^{k-1}g_{(k+1)(n-2)+2} $$
$$+\sum_{j=0}^{k-2}\frac{\varphi^j}{H}
g_{(k+1)(n-2)+k-j+3}+\varphi^{k-1}g_{(k+1)(n-2)+2}
+\varphi^{k}g_{(k+1)(n-2)+1} $$
$$+\sum_{j=0}^{k-1}\frac{\varphi^j}{H}
g_{(k+1)(n-2)+k-j+3}+\varphi^{k}g_{(k+1)(n-2)+1}
+\varphi^{k+1}g_{(k+1)(n-2)}. $$
An easy calculation yields that the above expression can be rewritten
in the form
$$q_{k+1}=\sum_{j=0}^{k}\frac{\varphi^j}{H^{k-j}}g_{(k+1)(n+1)-3j-2}
+\varphi^{k+1}g_{(k+1)(n-2)}.$$
Finally, it remains to use the fact we already established above
that $q_{k+1}^*$ (the coefficient at $dH$ in $\Omega^*_{k+1}$) is a
function of the same kind as the former $q_{k+1}$. $\Box$
%
%\vspace{2ex}
%\noindent
\begin{proposition}
\label{three}
 Assume that $M_1(t)=\ldots=M_k(t)\equiv 0$.
Then $\Omega_{k+1}=q_k\Omega_1=q_k\omega_n$ takes the form
$$\begin{array}{ll}
{\displaystyle \Omega_{k+1}=\alpha_{2n-2}(H)\sigma_0+\gamma_{2n-4}(H)\sigma_2
+\frac{a_0}{4H}(5\sigma_2-\sigma_0)+dQ_{k+1}+q_{k+1}dH} &
\mbox{\it if}\;\;k=1,\\[3mm]
{\displaystyle \Omega_{k+1}=\frac{\alpha_{(k+1)(n+1)-4}(H)}{H^{k-1}}\sigma_0+
\frac{\gamma_{(k+1)(n+1)-6}(H)}{H^{k-1}}\sigma_2 +dQ_{k+1}+q_{k+1}dH} &
\mbox{\it if}\;\;k>1.
\end{array}$$
\end{proposition}
%
%\vspace{2ex}
%\noindent
{\bf Proof.} We use formula (\ref{9}) from the proof of Proposition
\ref{two} and the fact that the function $\varphi$ is determined up
to an additive constant, say $c$. Recall that
$M_{k+1}(t)=\int_{\delta(t)}\Omega^*_{k+1}$ where $\Omega^*_{k+1}$ is
given by (\ref{9}). As above, one can use Lemma \ref{basic}
to express the last term in (\ref{9})
$$a_{k-1}\varphi^{k-1}xd\varphi=\frac{a_{k-1}}{H}\varphi^{k-1}\omega_3$$
as
$$\frac{a_{k-1}}{H}\{[\varphi^{k-1}(\alpha_2\sigma_0+\gamma_0\sigma_2)
+\mbox{\rm l.o.t} ] +dQ+qdH\}$$
where we denoted by l.o.t. the terms containing $\varphi^j$ with $j<k-1$.
The values of $\alpha_2$ and $\gamma_0$ can be calculated
from the last equation in (\ref{2}) which yields
\begin{equation}
\label{10}
a_{k-1}\varphi^{k-1}xd\varphi=\frac{a_{k-1}}{4H}\{[\varphi^{k-1}
(5\sigma_2-\sigma_0) +\mbox{\rm l.o.t} ]+dQ+qdH\}.
\end{equation}
Let us now put $\varphi+c$ instead of $\varphi$ in the formula of
$M_{k+1}(t)$. Then $M_{k+1}(t)$ becomes a polynomial in $c$ of degree $k$
with coefficients depending on $t$. Since $M_{k+1}$ does not depend on
this arbitrary constant $c$, all the coefficients at $c^j$, $1\leq j\leq k$
should vanish. By (\ref{9}), the coefficient at $c^k$ equals
$$\alpha_{(k+1)(n-2)+1}(t)I_0(t)+\gamma_{(k+1)(n-2)-1}(t)I_2(t)$$
which is zero as $M_{k+1}(t)$ does not depend on $c$.
This is equivalent to
$\alpha_{(k+1)(n-2)+1}(t)=\gamma_{(k+1)(n-2)-1}(t)\equiv 0.$
When $k=1$, this together with (\ref{9}) and (\ref{10}) implies the formula for
$\Omega_2$. Assume now that $k>1$.
When the leading coefficient at $c^k$ vanishes, the next coefficient,
at $c^{k-1}$, becomes
$$\left[\alpha_{(k+1)(n-2)+2}(t)-\frac{a_{k-1}}{4t}\right]I_0(t)
+\left[\gamma_{(k+1)(n-2)}(t)+5\frac{a_{k-1}}{4t}\right]I_2(t)$$
and both coefficients at $I_0$ and $I_2$ are identically zero which
yields $\alpha_{(k+1)(n-2)+2}=\gamma_{(k+1)(n-2)}\equiv 0$ and $a_{k-1}=0$.
Similarly, all coefficients in (\ref{9}) $\alpha_{(k+1)(n+1)-3j-4}$,
$\gamma_{(k+1)(n+1)-3j-6}$, $j>0$, become zero which proves
Proposition \ref{three}.  $\Box$

\vspace{2ex}
In the calculations above we took the eight-loop Hamiltonian
$H=\frac12y^2+\frac14(x^2-1)^2$
and considered the outer period annulus of the Hamiltonian vector field
$dH=0$, defined for levels $H=t$ with
$t\in\Sigma=(\frac14,\infty)$. Evidently a very minor modification (sign
changes in front of some terms in the formulas like (\ref{2})) is needed to
handle the double-heteroclinic Hamiltonian $H=\frac12y^2-\frac14(x^2-1)^2$
and the global-center Hamiltonian $H=\frac12y^2+\frac14(x^2+1)^2$.
The functions $\varphi$, $G$ and the interval $\Sigma$ could then be taken
respectively as follows:
$$\varphi=\frac{1}{2\sqrt2}\log\frac{1-x^2-\sqrt2y}{1-x^2+\sqrt2y},
\quad G=\frac{(x^2-1)y}{4},\quad {\textstyle\Sigma=(-\frac14,0)},$$
$$\varphi=\frac{i}{2\sqrt2}\log\frac{x^2+1+i\sqrt2y}{x^2+1-i\sqrt2y},
\quad G=\frac{(x^2+1)y}{4},\quad {\textstyle\Sigma=(\frac14,\infty)}.$$

Below, we state Theorem \ref{one} in a form to hold for all the three cases.
Recall that $I_k(t)=\int_{\delta(t)}\sigma_k= \int_{\delta(t)}x^kydx$,
$k=0,1,2$, where $\delta(t)$, $t\in\Sigma$, is the oval formed by the
level set $\{H=t\}$ for any of the three Hamiltonians.
%
%\vspace{2ex}
%\noindent
\begin{theorem}
\label{one}
For $t\in\Sigma$, the first nonvanishing generating
function $M_k(t)=\int_{H=t}\Omega_k$ corresponding to degree $n$
polynomial perturbations $dH-\varepsilon\omega_n=0$, has the form
$$\begin{array}{ll}
\mbox{\it for}\; k=1,& M_1(t)=\alpha_{\frac{n-1}{2}}(t)I_0(t)
+\gamma_{\frac{n-3}{2}}(t)I_2(t),\\[2mm]
\mbox{\it for}\; k=2,&{\displaystyle M_2(t)=
\frac{1}{t}\left[\alpha_n(t)I_0(t)
+\gamma_{n-1}(t)I_2(t)\right],}\\[2mm]
\mbox{\it for}\; k>2,&{\displaystyle M_k(t)=
\frac{1}{t^{k-2}}\left[\alpha_{\frac{k(n+1)}{2}-2}(t)I_0(t)
+\gamma_{\frac{k(n+1)}{2}-3}(t)I_2(t)\right]},
\end{array}$$
where $\alpha_j(t)$, $\gamma_j(t)$ denote polynomials in $t$ of
degree $[j]$.
\end{theorem}

\vspace{2ex}
\noindent
{\bf Proof.} Take a perturbation $dH-\varepsilon\omega_n=0$ where $\varepsilon$
is a small parameter. Then by a generalization of Fran\c{c}oise's recursive
procedure, one obtains $M_1(t)=\int_{\delta(t)}\Omega_1$,
and when $M_1(t)=\ldots=M_{k-1}(t)\equiv 0$, then
$M_k(t)=\int_{\delta(t)}\Omega_k$,  where
$\Omega_1=\omega_n$, $\Omega_k=q_{k-1}\Omega_1$ and $q_{k-1}$ is
determined from the representation $\Omega_{k-1}=dQ_{k-1}+q_{k-1}dH.$
The algorithm is effective provided we are able to express the
one-forms $\Omega_k$ in a suitable form which was done above.
For $k=1$, the result follows from (\ref{3}) applied with $m=n$.
For $k>1$, the result follows immediately from Proposition \ref{three}.
$\Box$

\vspace{2ex}
\noindent
Clearly, Theorem \ref{one} allows one to give an upper bound
to the number of zeros
of $M_k(t)$ in $\Sigma$ and thus to estimate from above the number of limit
cycles in the perturbed system which tend as $\varepsilon\to 0$  to periodic
orbits of the original system that correspond to Hamiltonian levels in
$\Sigma$. For this purpose, one can apply the known sharp results on
non-oscillation of elliptic integrals (most of them due to Petrov, see
also \cite{gi03}, \cite{rz91} and the references therein) to obtain the
needed bounds. Define the vector space
$${\cal M}_m=\{P_m(t)I_0(t)+P_{m-1}(t)I_2(t):\;
P_k\in\R[t],\;\mbox{\rm deg}\,P_k\leq k,\;t\in\Sigma\}.$$
Clearly, $\mbox{\rm dim}\,{\cal M}_m=2m+1$. We apply to the eight-loop case
Theorem 2.3 (c), (d) and Lemma 3.1  from \cite{rz91} and
to the double-heteroclinic and the global-center cases,
Theorem 2 (4), (5) and Lemma 1 (iii) from \cite{gi03} to obtain the
following statement.
%
%\vspace{2ex}
%\noindent
\begin{proposition}
\label{four}
{\em (i)} In the eight-loop case, any nonzero function
in ${\cal M}_m$ has at most $\mbox{\rm dim}\,{\cal M}_m=2m+1$
zeros in $\Sigma$.

\vspace{1ex}
\noindent
{\em (ii)} In the double-heteroclinic and the global-center cases, any
nonzero function in ${\cal M}_m$ has at most
$\mbox{\rm dim}\,{\cal M}_m-1=2m$ zeros in $\Sigma$.
\end{proposition}
%
%\vspace{1ex}
%\noindent
By Proposition \ref{four} and Theorem \ref{one}, we  obtain:
%
%\vspace{2ex}
%\noindent
\begin{theorem}
\label{t2}
 In the eight-loop case, the upper bound $N(n,k)$ to the
number of isolated zeros in $\Sigma$ of the first nonvanishing generating
function $M_k(t)$ corresponding to degree $n$ polynomial perturbations
$dH-\varepsilon\omega_n=0$, can be taken as follows:
$N(n,1)=2[\frac{n-1}{2}]+1$, $N(n,2)=2n+1$ and
$N(n,k)=2[\frac{k(n+1)}{2}]-3$ for $k>2$.
\end{theorem}
%
%\vspace{2ex}
%\noindent
\begin{theorem}
\label{t3}
 In the double-heteroclinic and the global-center cases,
the upper bound $N(n,k)$ to the number of isolated zeros in $\Sigma$
of the first nonvanishing generating
function $M_k(t)$ corresponding to degree $n$ polynomial perturbations
$dH-\varepsilon\omega_n=0$, can be taken as follows:
$N(n,1)=2[\frac{n-1}{2}]$, $N(n,2)=2n$ and
$N(n,k)=2[\frac{k(n+1)}{2}]-4$ for $k>2$.
\end{theorem}
Similarly, one can consider in the eight-loop case any of the internal
period annuli when the $(*)$ property holds. Take $t\in\Sigma=(0,\frac14)$
and consider the corresponding oval $\delta(t)$ lying (say) in the
half-plane $x>0$. Define the vector space
$$\textstyle
{\cal M}_m=\{P_{[\frac{m}{2}]}(t)I_0(t)+P_{[\frac{m-1}{2}]}(t)I_1(t)
+P_{[\frac{m-2}{2}]}(t)I_2(t):\;
P_k\in\R[t],\;\mbox{\rm deg}\,P_k\leq k,\;t\in\Sigma\}.$$
Clearly, $\mbox{\rm dim}\,{\cal M}_m=[\frac{3m+2}{2}]$.
By Petrov's result \cite{pet89}, any function in ${\cal M}_m$ has at most
$\mbox{\rm dim}\,{\cal M}_m-1=[\frac{3m}{2}]$ isolated zeros.
Applying this statement to (\ref{kth}), we get
%
%\vspace{2ex}
%\noindent
\begin{theorem}
\label{t4}
 In the internal eight-loop case, the number of isolated
zeros in $\Sigma$ of the first nonvanishing generating function $M_k(t)$
corresponding to degree $n$ polynomial perturbations
$dH-\varepsilon\omega_n=0$ is at most $N(n,k)=[\frac{3k(n-1)}{2}]$.
\end{theorem}
%
%\vspace{2ex}
%\noindent
It is well known that the bounds in Theorems
\ref{t2}, \ref{t3}, \ref{t4} are sharp for $k=1$.
That is, there are degree $n$ perturbations with the prescribed numbers
of zeros of $M_1(t)$ in the respective $\Sigma$. One cannot expect that
this would be the case for all $k>1$ and $n$. The reason is that $M_k$,
$k>1$, is a very specific function belonging to the linear space ${\cal M}_m$
with the respective index $m$ which in general would not possess the maximal
number of zeros allowed in ${\cal M}_m$. Moreover, as there is a finite
number of parameters in any $n$-th degree polynomial perturbation,
after a finite steps the perturbation will become an integrable one
and hence $M_k(t)$ will be zero for all $k>K$ with a certain (unknown)
$K$. The determination of the corresponding $K$ and the exact upper bound
to the number of isolated zeros that the functions from the set
$\{M_k(t):\; 1\leq k\leq K\}$ can actually have in $\Sigma$, are huge
problems. We will not even try to solve them here. Instead, below we show
that the result in Theorem \ref{one} can be slightly improved when $k>1$ and
$n$ is odd.

\vspace{2ex}
\noindent
{\bf Theorem 3$^+$} {\it For $t\in\Sigma$ and $n$ odd, the first nonvanishing
generating function $M_k(t)=\int_{H=t}\Omega_k$ corresponding to degree $n$
polynomial perturbations $dH-\varepsilon\omega_n=0$, has the form
$$\begin{array}{ll}
\mbox{\it for}\; k=1,& M_1(t)=\alpha_{\frac{n-1}{2}}(t)I_0(t)
+\gamma_{\frac{n-3}{2}}(t)I_2(t),\\[2mm]
\mbox{\it for}\; k=2,&{\displaystyle M_2(t)=
\frac{1}{t}\left[\alpha_{n-1}(t)I_0(t)
+\gamma_{n-1}(t)I_2(t)\right],}\\[2mm]
\mbox{\it for}\; k>2,&{\displaystyle M_k(t)=
\frac{1}{t^{k-2}}\left[\alpha_{\frac{k(n+1)}{2}-3}(t)I_0(t)
+\gamma_{\frac{k(n+1)}{2}-4}(t)I_2(t)\right]},
\end{array}$$
where $\alpha_j(t)$, $\gamma_j(t)$ denote polynomials in $t$ of degree $j$.}

\vspace{2ex}
\noindent
{\bf Proof.} Given $A(x,y,H)$, a polynomial of weighted degree $m$, we
denote by $\bar{A}$ its highest-degree part:
$$\bar{A}(x,y,H)=\sum_{i+j+2k=m}a_{ijk}x^iy^jH^k.$$
The same notation will be used for the respective polynomial one-forms.
We begin by noticing that
$$\bar\omega_n=(a_0y^n+a_1xy^{n-1}+a_2x^2y^{n-2})dx+d(b_0y^{n+1}+b_1xy^n
+b_2x^2y^{n-1}+b_3x^3y^{n-2})$$
because all terms containing $x^j$ with $j\geq 4$ can be expressed through
lower-degree terms. If $M_1(t)\equiv 0$ then, by Proposition \ref{petrov},
$\bar\alpha_{n-1}=\tilde\gamma_{n-3}=0$ which implies that $a_0=a_2=0$,
see equations (\ref{suppl}). From the formulas we derived in the proof
of Proposition \ref{petrov}, one can also obtain that, up to a
lower-degree terms,
$$xy^{n-1}dx=\frac{2(n-1)}{n}Hxy^{n-3}dx
-\frac{n-1}{2n}x^2y^{n-3}dH+d\frac{x^2y^{n-1}}{2n},$$
which yields
$$xy^{n-1}dx=dx^2P_{n-1}(y,H)-x^2P_{n-3}(y,H)dH+ \mbox{\rm l.d.t.}$$
where $P_j$ denotes a weighted homogeneous polynomial of weighted degree
$j$ with positive coefficients. Now,
$$\bar{\Omega}_2=\overline{q_1\omega}_n=-a_1x^2P_{n-3}d(b_0y^{n+1}+b_1xy^n)$$
and we see that the highest-degree coefficient of the polynomial $\alpha_n(t)$
in the formula of $M_2(t)$ should be zero.
If, in addition, $M_2(t)\equiv 0$, then $a_1b_1=0$. When $a_1=0$, one obtains
$\bar{q}_1=0\Rightarrow\bar\Omega_k=0$, $k\geq 2$ and the claim follows.
If $b_1=0$, then
$\bar{\Omega}_2$ is proportional to $x^2P_{n-3}y^ndy$ which implies that
all $\bar{q}_k$, $k\geq 2$, will have the form
$\bar{q}_k=x^2P_{k(n-1)-2}(y,H)$ where $P_j$ are as above, and hence,
$\bar{\Omega}_{k+1}=\overline{q_k\omega}_n$ will have no impact on the
value of $M_{k+1}$. $\Box$

\vspace{1ex}
The result in Theorem 3$^+$ allows one to improve Theorems \ref{t2} and
\ref{t3}, but we are not going to present here the obvious new statements.

\subsection{The $D_4$ singularity}
\begin{figure}
\begin{center}
\epsfig{file=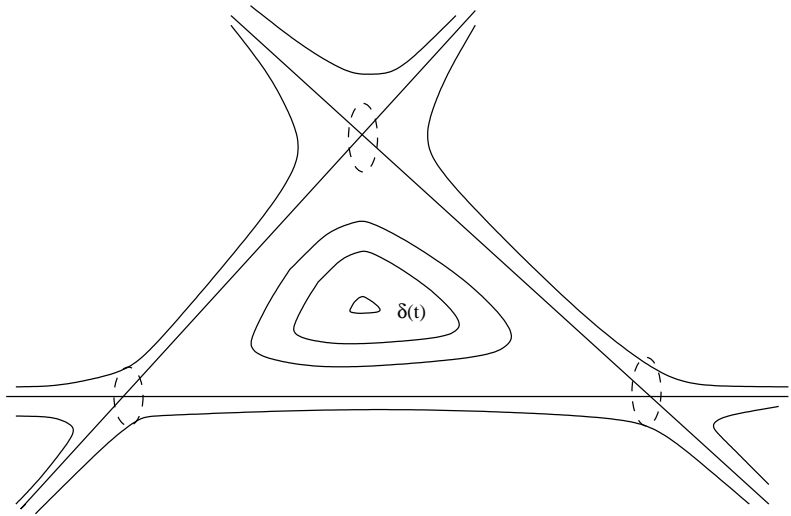}
\end{center}
\caption{ The level sets of $f= x[y^2-(x-3)^2]$ and the family of ovals
$\delta (t)$.}
\label{triangle}
\end{figure}
Let
$$
f= x[y^2-(x-3)^2]
$$
and denote by $\delta (t)$ the family of ovals defined by $\{(x,y)\in
\R^2: f(x,y)=t\}$,
$t\in (-4,0)$, see Fig. \ref{triangle}.
 We will denote by the same letters the
corresponding continuous families of free homotopy classes of loops defined
on the universal covering space of $\C \setminus \{0,-4\}$,
and fix $t_0 \neq 0,-4$.
\begin{proposition}
\label{pr2}
We have
$$
H^{\delta(t_0)}_1(f^{-1}(t_0),\Z)=\Z^3
$$
and the kernel of the canonical map $H^{\delta
_e(t_0)}_1(f^{-1}(t_0),\Z) \rightarrow H_1(f^{-1}(t_0),\Z)$ is
equal to $\Z$.
\end{proposition}
{\bf Proof.}
\begin{figure}
\begin{center}
\epsfig{file=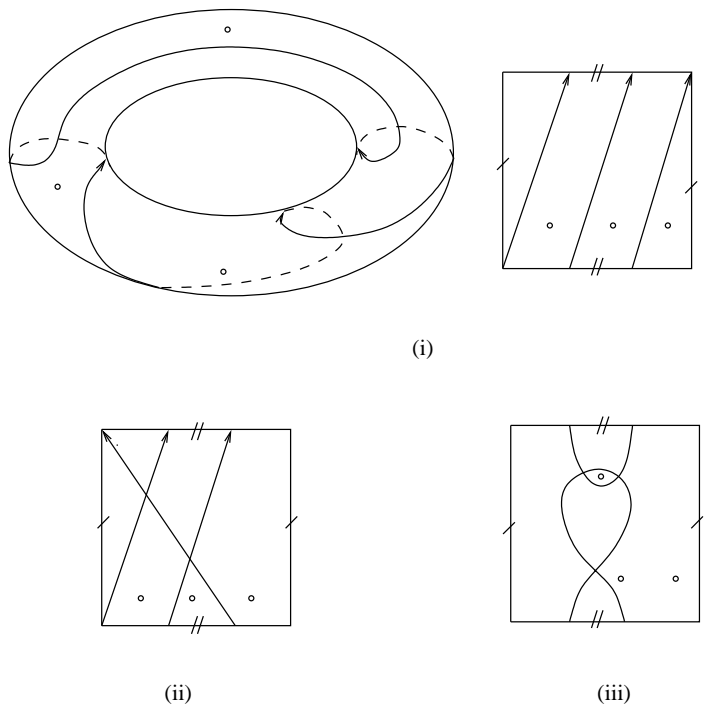}
\end{center}
\caption{ (i) The closed loop
 $l_{0*} \delta (t_0)$, (ii)
$\mbox{\it Var}_{l_{0}}\delta (t_0)$,
and (iii) $\mbox{\it Var}^2_{l_{0}}\delta (t_0)$.}
\label{g3}
\end{figure}
The fibers $f^{-1}(t) \subset \C^2$ for $t\neq 0,-4$ are genus-one surfaces
with three removed points. Let
$l_0,l_{-4} \in \pi _1(\C\setminus \{0,-4\},t_0)$ be two simple loops
making one turn around $0$
and $-4$ respectively in a positive direction. The closed loop
 $l_{0*} \delta (t_0)$ is shown on Fig. \ref{g3}, (i). The loops representing
$\mbox{\it Var}_{l_{0}}\delta (t_0)$, $\mbox{\it
Var}^2_{l_{0}}\delta (t_0)$, where $\mbox{\it Var}_{l_{0}}=
(l_{0}-id)_*$, are shown on Fig. \ref{g3}, (ii), (iii)
respectively. It follows that $\mbox{\it Var}^3_{l_{0}}\delta
(t_0)$ may be represented by a loop homotopic to a point. Finally,
the variation of an arbitrary element of $H_1^\delta(f^{-1}(t_0)_t
,\Z) $ along $l_{-4}$ is a composition of free homotopy classes of
$\delta $ (several  times) which shows that
$H_1^\delta(f^{-1}(t_0)_t ,\Z) $ is generated by
$$
\delta(t_0) ,\;\; \mbox{\it Var}_{l_{0}}\delta (t_0), \;\;
\mbox{\it Var}^2_{l_{0}}\delta (t_0) \; .
$$
The equivalence class $\mbox{\it Var}^2_{l_{0}}\delta (t_0)$ is
homologous to zero while the other two are homologically
independent. This shows that the image of
$H_1^\delta(f^{-1}(t_0)_t ,\Z) $ in $H_1(f^{-1}(t_0),\Z)$ is
$\Z^2$. It remains to show that the equivalence class of $k
\mbox{\it Var}^2_{l_{0}}\delta (t_0)$ in $H_1^\delta(f^{-1}(t_0)_t
,\Z) $ is nonzero for any $k\in \Z$.

\begin{figure}
\begin{center}
\epsfig{file=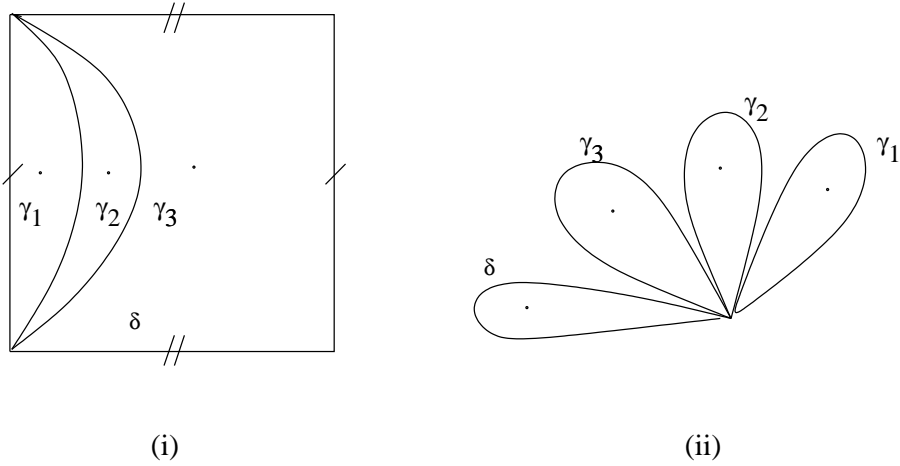}
\end{center}
\caption{ (i) The generators of the fundamental group $ \pi
_1(f^{-1}(t_0),P_0)$ ; (ii)
The generators of the fundamental group $ \pi _1(\C \setminus
\{z_0,z_1,z_2,z_3 \},\tilde{z})$ .
}
\label{generators}
\end{figure}
The fundamental group $ \pi _1(f^{-1}(t_0),P_0)$ is a free group
with generators $\delta, \gamma _1, \gamma _2, \gamma_3$ shown on
Fig. \ref{generators}, (i). We have
\begin{equation}\label{var}
\mbox{\it Var}_{l_{0}}\delta (t_0)= \gamma _1 \gamma _2 \gamma_3,
\;\; \mbox{\it Var}^2_{l_{0}}\delta (t_0) = \gamma _1\gamma
_2\gamma _1^{-1}\gamma _2^{-1} .
\end{equation}
Let
$$
S= \{\delta ,\gamma _1 \gamma _2 \gamma_3,
[\gamma _1,\gamma _2]\} \mbox{  where  } [\gamma _1,\gamma _2]= \gamma
_1\gamma _2\gamma _1^{-1}\gamma _2^{-1}\}
$$
and let $\hat{S}$ be the least normal subgroup of $ \pi
_1(f^{-1}(t_0),P_0)$ containing $S$. A general method to study
$H_S= \hat{S}/[\hat{S},\pi _1(f^{-1}(t_0),P_0)]$ consists of
constructing its dual space. Namely, let $z_0,z_1,z_2,z_3$ be
distinct complex numbers and let $\delta ,\gamma _1,\gamma
_2,\gamma_3$ be simple loops making one turn about
$z_0,z_1,z_2,z_3$ respectively in a positive direction as it is
shown on Fig. \ref{generators}, (ii). Note that
$$\pi _1(\C\setminus \{z_0,z_1,z_2,z_3\}, \tilde{z})=
\pi _1(f^{-1}(t_0),P_0) \; .
$$
Let
$$
\omega = \ln \frac{z-z_1}{z-z_3}\left(\frac{1}{z-z_2} - \frac{1}{z-z_1}\right)
dz \; .
$$
We claim that $\omega $ defines a linear function on $H_S$ by the formula
$$
l \rightarrow \int_l \omega \; .
$$
Indeed,  whatever  the determination of the multivalued function {\it ln} be,
we have $\int_{\delta }\omega  = 0$, and
$\int_{\gamma _1\gamma _2\gamma _3 }\omega  $ is well
defined. The latter holds true because
$$
\int_{\gamma _1\gamma _2\gamma _3 }\left(
\frac{1}{z-z_2} - \frac{1}{z-z_1}\right) dz = 0
$$
and $\ln \frac{z-z_1}{z-z_3}$ is single-valued along the loop
$\gamma _1\gamma _2\gamma _3 $. Finally, along $[\gamma _1,\gamma _2]$
the differential $\omega $ is single-valued too and
$\int_{[\gamma _1,\gamma _2]} \omega $ does not depend on the
determination of $\omega $. An easy exercise
shows that $\int_{[\gamma _1,\gamma _2]} \omega = -4\pi ^2$. We
conclude that the space dual to $H_S$ is generated (for instance)
by $\omega , dz/(z-z_0), dz/(z-z_1)$ and hence $H_S=\Z^3$.
Obviously the kernel of the homomorphism $H^{\delta
_e(t_0)}_1(f^{-1}(t_0),\Z) \rightarrow H_1(f^{-1}(t_0),\Z)$ is the
infinite cyclic group generated by the commutator $[\gamma
_1,\gamma _2]$. $\square$

\vspace{2ex}
\noindent
According to Theorem \ref{main} and Proposition \ref{pr2} the generating
function $M(t)$ might not be an Abelian integral, the obstruction being
 the kernel of the map
$H^{\delta _e(t_0)}_1(f^{-1}(t_0),\Z) \rightarrow
H_1(f^{-1}(t_0),\Z)$.
 Indeed, it follows from \cite{ili96}, \cite{zol94} that for some quadratic
unfoldings of $\{df=0\}$, the corresponding generating function
$M_{\delta (t)}$ is not an Abelian integral (see the open question 3. at
the end of section
\ref{mainresult}).
More explicitly, we have
\begin{proposition}
The generating function associated to the unfolding
$$\textstyle df+\varepsilon (2-x+\frac12x^2)dy=0,\quad f=x[y^2-(x-3)^2]$$
and to the family of ovals around the center of the unperturbed
system, is not an Abelian integral of the form $(\ref{ai})$.
It satisfies an equation of Fuchs type of order three.
\end{proposition}
{\bf Proof}. For a convenience of the reader, below we
present the needed calculation. Denote $\omega_2=-(2-x+\frac12x^2)dy$.
One can verify  \cite{advde} that $\omega_2=dQ_1+q_1df$, with
$${\textstyle Q_1=\frac16[fL(x,y)-x^2y-12y],\quad
q_1=-\frac16L(x,y),} \quad L(x,y)=\ln\frac{3-x-y}{3-x+y},$$
and that the form $q_1\omega_2-q_2df$ is exact, where
$$q_2=\frac{L^2}{72}+\frac{x^3-3x^2+12x-36}{36f}$$
(to check this, we make use of the identity $fdL=2xydx+(6x-2x^2)dy$).
Therefore $M_1(t)=M_2(t)\equiv 0$ for this perturbation, and
$$M_3(t)=\int_{\delta(t)}q_2\omega_2=\int_{\delta(t)}q_2dQ_1=
\frac{1}{216}\int_{\delta(t)}(x^3-3x^2+12x-36)dL$$
$$+\frac{1}{216}\int_{\delta(t)}(x^2+12)yd\left(\frac{L^2}{2}
+\frac{x^3-3x^2+12x-36}{t}\right).$$
In the same way as in \cite{ili96}, Appendix, we then obtain
$$M_3(t)=\frac{1}{36t}\int_{\delta(t)}{\textstyle
[36(x-1)\ln x +\frac12x^4-\frac72x^3-\frac{39}{2}x^2+12x+24]ydx.}$$
As $I_1=I_0$ and $(2k+6)I_{k+1}=(12k+18)I_k-18kI_{k-1}-(2k-3)tI_{k-2}$,
the final formula becomes
$$M_3(t)=\frac{1}{t}\int_{\delta(t)}y(x-1)\ln x dx
-\frac{3}{32}\int_{\delta(t)}\frac{ydx}{x}.$$
For a general quadratic perturbation satisfying $M_1(t)=M_2(t)\equiv 0$,
the formula of $M_3(t)$ will take the form \cite{ili96}, \cite{advde}
\begin{equation}
\label{M3}
M_3(t)=c_{-1}I_{-1}(t)+\left(c_0+\frac{c_1}{t}\right)I_0
+\frac{c_*}{t}I_*(t),\quad I_*(t)=\int_{\delta(t)}y(x-1)\ln x dx
\end{equation}
where $c_j, c_*$ are some constants depending on the perturbation.
Below we write up the equation satisfied by $M_3(t)$ and show that,
apart of $M_1$ and $M_2$, $M_3$ is not an Abelian integral,
due to $I_*$. We can rewrite (\ref{M3}) as
$tM_3(t)=(\alpha+\beta t)I_0+\gamma I_2+\delta I_*$ (with some
appropriate constants) and
use the Fuchsian system satisfied by ${\bf I}=(I_*, I_2, I_0)^\top$
\cite{ili96}, namely
$${\bf I}={\bf A}{\bf I}',\quad\mbox{where}\quad
{\bf A}=\pmatrix{t & -2 & t+6\cr
0 & \frac34(t-6) & \frac32(t+9)\cr
0 & -3 & \frac32(t+6)}$$
to derive explicitly the third-order Fuchsian equation satisfied by $M_3(t)$.
One obtains
$$DP(t^2M_3')''+(tP-DP')(t^2M_3')'+Q(t^2M_3')=0,$$
where  $D=t(t+4)$ and
$$\begin{array}{rl}
P&\!\!=(8\beta^2-\beta\gamma)t^3
-(56\alpha\beta+\alpha\gamma+96\beta\gamma+2\gamma^2
+48\beta\delta+2\gamma\delta)t^2\\
&\!\!+(8\alpha^2-288\alpha\beta+12\alpha\gamma-432\beta\gamma+24\alpha\delta
-192\beta\delta+28\gamma\delta+16\delta^2)t\\
&\!\!+(96\alpha\delta+144\gamma\delta+64\delta^2),\\
Q&\!\!=\frac49\{(40\beta^2-5\beta\gamma)t^3-(64\alpha\beta+2\alpha\gamma
-288\beta^2+144\beta\gamma+4\gamma^2 +48\beta\delta+4\gamma\delta)t^2\\
&\!\!+(4\alpha^2-144\alpha\beta+12\alpha\gamma-432\beta\gamma+12\alpha\delta
-240\beta\delta-4\gamma\delta+8\delta^2)t+32\delta^2\}.\end{array}$$
 For the above particular perturbation, the equation of $M_3$ reads
% $$t(t+4)(39t^2+704t+2048)(t^2M_3')''-t(39t^2+312t+768)(t^2M_3')'$$
% $$+{\textstyle\frac89}(39t^2+311t+512)(t^2M_3')=0$$
% or equivalently,
$$
t^2(t+4)(39t^2+704t+2048)M_3'''+t(117t^3+3128t^2+18688t+32768)M_3''
$$
$$
+{\textstyle\frac89}(39t^3+1544t^2+9728t+18432)M_3' = 0.\;\;
$$
The above equation is obviously of Fuchs type and its monodromy
group is studied in a standard way. The characteristic exponents
associated to the regular singular point $t=0$ are $-1,0,0$.
Further analysis (omitted) shows that the monodromy transformation
of a suitable fundamental set of solutions along a small closed
loop about $t=0$ reads
$$
 \pmatrix{
1 & 1& 0  \cr
0 & 1& 1 \cr
0 & 0& 1  } .
$$
Indeed, according to formula (19) in \cite{ili96} in a
neighborhood of $t=0$ we have
$$
I_*(t)=\int_{\delta(t)}y(x-1)\ln x dx= -6 - \frac{1}{6} t \ln^2 t
+ \dots
$$
From this  we obtain that $\mbox{\it Var}^2_{l_{0}}M_3(t)
\not\equiv 0. $ On the other hand
$$
\mbox{\it Var}_{l_0}^2 M_3(\delta(t_0), {\cal F}_\varepsilon ,t)=
M_3(\mbox{\it Var}_{l_0}^2 \delta(t_0), {\cal F}_\varepsilon ,t)
$$
where the loop $\mbox{\it Var}_{l_0}^2 \delta(t_0)= \gamma
_1\gamma _2\gamma _1^{-1}\gamma _2^{-1} $ is homologous to zero,
see (\ref{var}). If $M_3$ were an Abelian integral then its second
``variation" $M_3(\mbox{\it Var}_{l_0}^2 \delta(t_0), {\cal
F}_\varepsilon ,t)$ would vanish identically which is a
contradiction.
 $\Box$

\vspace{2ex}
\noindent
%%%%%%%%%%%%%
{\bf Acknowledgments.}
We are grateful to the referee for the useful remarks and valuable
recommendations. The second author has been supported by the Ministry
of Science and Education of Bulgaria, Grant no. MM1403/2004.
%%%%%%%%%%%%%

\end{document}